# ISOPERIMETRIC INEQUALITIES AND MIXING TIME FOR A RANDOM WALK ON A RANDOM POINT PROCESS

BY PIETRO CAPUTO AND ALESSANDRA FAGGIONATO

*Università di Roma Tre and Università di Roma "La Sapienza"*

We consider the random walk on a simple point process on $\mathbb{R}^d$, $d \geq 2$, whose jump rates decay exponentially in the $\alpha$-power of jump length. The case $\alpha = 1$ corresponds to the phonon-induced variable-range hopping in disordered solids in the regime of strong Anderson localization. Under mild assumptions on the point process, we show, for $\alpha \in (0, d)$, that the random walk confined to a cubic box of side $L$ has a.s. Cheeger constant of order at least $L^{-1}$ and mixing time of order $L^2$. For the Poisson point process, we prove that at $\alpha = d$, there is a transition from diffusive to subdiffusive behavior of the mixing time.

**1. Introduction.** We consider the following model of a random walk in random environment. We let $\xi$ denote the realization of a simple point process on $\mathbb{R}^d$, $d \geq 2$, and identify $\xi$ with the countable collection of its points. We study the continuous-time Markov chain with state space $\xi$ and with jump rate from $x$ to $y$ given by a negative exponential of the Euclidean distance to some power $\alpha > 0$:

$$(1.1) \qquad r_{x,y} = e^{-|x-y|^\alpha}, \qquad x \neq y \text{ in } \xi.$$

The canonical examples are obtained when $\xi$ is the realization of a homogeneous Poisson point process, but our assumptions on the environment will allow more general processes. We consider the random walk obtained by confining the particle to a cubic box $\Lambda_L \subset \mathbb{R}^d$ with side $L$, that is, the random walk with rates (1.1) on $\xi_L := \xi \cap \Lambda_L$. There are two natural processes associated with the rates (1.1), defined as follows. Set $w_x := \sum_{y \in \xi_L \setminus \{x\}} r_{x,y}$, $x \in \xi_L$. In the first model, the particle at $x$ waits an exponential time with mean $1/w_x$ and then jumps to some $y \in \xi_L \setminus \{x\}$ with probability $r_{x,y}/w_x$.

In the second model, the particle at $x$ waits an exponential time with mean 1 and then jumps to $y \in \xi_L \setminus \{x\}$ with probability $r_{x,y}/w_x$. While the first model has a uniform stationary measure, in the second model, the stationary measure is given by the weights $w_x$ (see below). Our main results will concern both models, showing that they essentially share the same features.

If $\xi$ is the regular grid $\mathbb{Z}^d$, then it follows from well-known facts that the Cheeger constant or *conductance* of the Markov chains defined above is at least of order $L^{-1}$ and that the *mixing time* is of order $L^2$ (see, e.g., [25]). We shall show that, for dimension $d \geq 2$, $\alpha < d$, these diffusive-type estimates continue to hold in our setting for typical realizations of the underlying point process. We establish our results by way of estimates on the so-called *isoperimetric profile* of the random walk. This will be achieved by combining stochastic domination and percolation techniques. Similar results have recently been obtained in the case of a random walk on the supercritical percolation cluster [2, 20]. The case $\alpha > d$, with $\xi$ a Poisson point process, will be shown to be subdiffusive in the sense that Cheeger's constant is smaller than any inverse power of $L$, thus implying that the mixing time is larger than any power of $L$. In the critical case $\alpha = d$, we find a transition from subdiffusive to diffusive behavior of the mixing time as the intensity of the process increases. In a separate work [6], we analyze the one-dimensional case in detail. Before describing our results, let us add a few more lines to motivate our work.

In recent work [10], a variant of the first model with $\alpha = 1$ has been studied, where each point $x \in \xi$ is given an independent random energy $E_x$ and the rates in (1.1) are multiplied by a Boltzmann-type factor involving the initial and final energies $E_x, E_y$, at inverse temperature $\beta$. This can be seen as a model for the study of the phonon-induced hopping conductivity observed in disordered solids in the regime of strong Anderson localization. Points $x \in \xi$ correspond to impurities of the solid and the quantum electron Hamiltonian has exponentially localized eigenfunctions with localization centers $x$ if the corresponding energy $E_x$ is near the Fermi level. The DC conductivity of such materials would vanish if it were not for lattice vibrations (phonons) at nonzero temperature. These induce transitions between the localized eigenstates, whose probability rate can be derived from the Fermi golden rule. Due to the localization, one can think of electrons as classical particles. Moreover, at large $\beta$ and within a mean field approximation, the resulting stochastic hopping dynamics is given by the random walk mentioned above. We refer to [10] and references therein for a thorough discussion of the physics behind the model.

Under suitable assumptions on the underlying point process, the authors of [10] obtain an invariance principle for the random walk and prove, in dimension $d \geq 2$, a lower bound on the effective diffusivity which coincides with the prediction of Mott law, that is, a power-law behavior of the logarithm of conductivity as a function of $\beta$. An upper bound in agreement with Mott law is proved in [9].

The invariance principle of [10] is based on classical homogenization results [8] which allow one to prove that the law of the rescaled random walk converges to the law of a Brownian motion, *in probability*, with respect to the environment. To prove *almost sure* convergence, a different approach is required; see [4, 18] for two different ways of obtaining the almost sure invariance principle in the case of a random walk on the supercritical percolation cluster. The situation in our continuum model is slightly different and we shall come back to the almost sure convergence problem in future work. One of the byproducts of the isoperimetric inequalities we establish in the present paper is the almost sure Poincaré inequality for finite boxes which, according to the approach of [18], may be seen as a first step in the program of proving an almost sure invariance principle. In this context, the introduction of the energy marks $E_x$ at the points $x \in \xi$ does not cause any additional technical difficulty and we will restrict our attention to the model defined by (1.1) corresponding to the case of all energies $E_x$ being equal. We remark that, in contrast with the one-dimensional case [6], for $d \geq 2$, an invariance principle with positive effective diffusion constant may well hold, even in the absence of diffusive bounds for the Poincaré constant; see [5, 19] for recent interesting developments on $\mathbb{Z}^d$ random walks among random conductances.

Finally, we point out some technical features of the random walks considered here that make their analysis somewhat subtle, especially for a fixed environment: The random walk is genuinely nonuniformly elliptic, the particle can perform arbitrarily large jumps and, when visiting a very isolated region of $\xi$ (with at least two points in the second model), spends much time there, but, when leaving such a region, performs a very long jump. This trapping effect will become particularly clear in the analysis of the random walk on a Poisson point process, where the transition from diffusive to subdiffusive behavior of the mixing time comes from trapping in isolated regions.

1.1. *Main results.* The main assumptions on the point process are as follows. We consider a simple point process $\mathcal{P}$ on $\mathbb{R}^d$, $d \geq 2$, that is, a probability measure $\mathcal{P}$ on the set $\Omega$ of locally finite subsets $\xi$ of $\mathbb{R}^d$, endowed with the $\sigma$-algebra $\mathcal{F}$ generated by the counting variables $N_\Lambda : \xi \to \#(\xi \cap \Lambda)$ (cardinality function), $\Lambda$ a bounded Borel subset of $\mathbb{R}^d$. We refer to [7] for a basic reference on point processes.

We write $\mathcal{P}_{*,\rho}$ if $\mathcal{P}$ is the homogeneous Poisson point process with intensity $\rho > 0$. Given a realization $\xi$ of the point process and a bounded Borel subset $\Lambda \subset \mathbb{R}^d$, we shall often write $\xi(\Lambda) = N_\Lambda(\xi)$ for the number of points of $\xi$ belonging to $\Lambda$. For any $K \in \mathbb{R}_+$, we write $Q_K = [0, K)^d$ for the cube of side $K$ in $\mathbb{R}^d$. Consider the partition of $\mathbb{R}^d$ into translates of $Q_K$, that is, $\mathbb{R}^d = \bigcup_{x \in \mathbb{Z}^d} B_x$, $B_x := xK + Q_K$, and declare a box $B_x$ *good* iff $\xi(B_x) \geq 1$. The configuration of good boxes may be described by the random field $\sigma = \sigma(K) := \{\sigma_x(\xi)\}_{x \in \mathbb{Z}^d}$ defined by

(1.2) $$\sigma_x(\xi) := \chi(\xi(B_x) \geq 1),$$

where $\chi(A)$ denotes the indicator function of the event $A$. We can now state our main assumption.

ASSUMPTION (A1). We say that $\mathcal{P}$ satisfies Assumption (A1) if, for every $0 < p < 1$, there exists $K \in \mathbb{R}_+$ such that the random field $\sigma(K)$ defined in (1.2) stochastically dominates the independent Bernoulli process on $\mathbb{Z}^d$ with parameter $p$.

We recall that the above statement is equivalent to saying that we can construct the process $\sigma = \sigma(K)$ and the independent Bernoulli process $Z$ on $\mathbb{Z}^d$ with parameter $p$ on the same probability space in such a way that $\sigma_x \geq Z_x$ almost surely. We refer to [15] for more on stochastic domination by Bernoulli product measures.

It is easy to check that $\mathcal{P}_{*,\rho}$ satisfies Assumption (A1) for any $\rho > 0$. We shall see (in Section 5) that this assumption is also satisfied by processes with nontrivial correlation structure, provided a suitable mixing condition is satisfied. We will also consider the stochastic order between point processes defined in the standard way (see, e.g., [12]). For two processes $\mathcal{P}, \mathcal{P}'$, we write $\mathcal{P}' \preceq \mathcal{P}$ if $\mathcal{P}'$ is stochastically dominated by $\mathcal{P}$. This is equivalent to the existence of a coupling of the fields $\mathcal{P}, \mathcal{P}'$ such that, almost surely, $\xi' \subset \xi$, with $(\xi, \xi')$ denoting the random sets with marginal distributions given by $\mathcal{P}$ and $\mathcal{P}'$, respectively. In particular, it follows that if there exists $\rho > 0$ such that $\mathcal{P} \succeq \mathcal{P}_{*,\rho}$, then $\mathcal{P}$ satisfies Assumption (A1).

For every $L \in \mathbb{N}$, $\Lambda_L$ is the cubic box centered at the origin and $\xi_L$ denotes the restriction of the process to $\Lambda_L$, that is,

$$\Lambda_L = \left[-\frac{L}{2}, \frac{L}{2}\right]^d, \qquad \xi_L = \xi \cap \Lambda_L.$$

Before stating the results, we need another regularity assumption which guarantees that local fluctuations in the number of points are not too large in our process. For every bounded Borel set $A \subset \mathbb{R}^d$, define

(1.3) $$R_A(\xi) = \sum_{x \in \xi \cap A} \sum_{y \in \xi} e^{-|x-y|^\alpha}.$$

We shall consider, for $\varepsilon \in (0,1)$, the cubes $V_x := L^\varepsilon x + [0, L^\varepsilon)^d$, $x \in \mathbb{Z}^d$, with side $L^\varepsilon$. We write $J_L$ for the set of $x \in \mathbb{Z}^d$ such that $V_x \cap \Lambda_L \neq \varnothing$.

ASSUMPTION (A2). We say that $\mathcal{P}$ satisfies Assumption (A2) if, for all $\varepsilon \in (0,1)$, there exists $C < \infty$ such that $\mathcal{P}$-almost surely

$$R_{V_x} \leq C L^{\varepsilon d} \qquad \forall x \in J_L, \tag{1.4}$$

for $L$ sufficiently large.

Here, and in all our statements below, we use the following convention: Given a sequence of events $\{E_L\}_{L \in \mathbb{N}}$, we say that $E_L$ holds $\mathcal{P}$-a.s. for $L$ sufficiently large if there exists $\Omega_0$ with $\mathcal{P}(\Omega_0) = 1$ such that for every $\xi \in \Omega_0$, there is $L_0(\xi) < \infty$ such that $\xi \in E_L$ for all $L \in \mathbb{N}$ with $L \geq L_0(\xi)$. In particular, the statement (1.4) says that if we define the event

$$E_L = \{\xi : R_{V_x}(\xi) \leq C L^{\varepsilon d}, \forall x \in J_L\},$$

then $E_L$ eventually holds $\mathcal{P}$-almost surely. Assumption (A2) is satisfied by Poisson point processes (see Section 5) and therefore by any process $\mathcal{P}$ such that $\mathcal{P} \preceq \mathcal{P}_{*,\rho}$ for some finite intensity $\rho$.

Let us remark that a consequence of Assumptions (A1) and (A2) is that there exists $C < \infty$ such that the inequalities

$$C^{-1} L^d \leq \xi(\Lambda_L) \leq C L^d \tag{1.5}$$

hold $\mathcal{P}$-a.s. for $L$ sufficiently large, where $\xi(\Lambda_L)$ is the number of points of $\xi$ in $\Lambda_L$. Indeed, the lower bound is a simple consequence of Assumption (A1), while the upper bound follows immediately from Assumption (A2) and the obvious fact that $R_A(\xi) \geq \xi(A)$ for any bounded Borel set $A \subset \mathbb{R}^d$.

We now define the random walk models to be considered. In addition to the two natural models previously discussed (which correspond, resp., to cases $i = 1$ and $i = 2$ below), we find it useful to introduce a third model. For $i = 1, 2, 3$ and $x, y \in \xi_L$ with $x \neq y$, we define

$$\mathcal{L}^i(x,y) = \frac{e^{-|x-y|^\alpha}}{w_x^i},$$

$$w_x^i := \begin{cases} 1, & i = 1, \\ \displaystyle\sum_{z \in \xi_L : z \neq x} e^{-|x-z|^\alpha}, & i = 2, \\ \max\left\{1, \displaystyle\sum_{z \in \xi_L : z \neq x} e^{-|x-z|^\alpha}\right\}, & i = 3. \end{cases} \tag{1.6}$$

On the diagonal, we set $\mathcal{L}^i(x,x) = -\sum_{z \in \xi_L : z \neq x} \mathcal{L}^i(x,z)$. For each $i = 1, 2, 3$, we write $X_t^i = X_t^i(\xi_L)$ for the continuous-time Markov chain with state space $\xi_L$ and infinitesimal generator $\mathcal{L}^i(x,y), x, y \in \xi_L$. Setting

$$\nu^i(x) := \frac{w_x^i}{\sum_{z \in \xi_L} w_z^i}, \qquad x \in \xi_L, \tag{1.7}$$

we obtain a symmetrizing probability for the matrix $\mathcal{L}^i$ and therefore the reversible probability measure for $X_t^i$ is given by $\nu^i$. Note that $\nu^1$ is uniform: $\nu^1(x) = \frac{1}{\xi(\Lambda_L)}$.

For any nonempty subset $U \subset \xi_L$, we define $U^c := \xi_L \setminus U$. The constant $I_U^i$, $i = 1, 2, 3$, is given by

$$I_U^i := \frac{\sum_{x \in U, y \in U^c} \nu^i(x) \mathcal{L}^i(x, y)}{\sum_{x \in U} \nu^i(x)}. \tag{1.8}$$

With the notation

$$W^i(U) := \sum_{x \in U} w_x^i, \tag{1.9}$$

we obtain

$$I_U^i = \frac{\sum_{x \in U, y \in U^c} e^{-|x-y|^\alpha}}{W^i(U)}. \tag{1.10}$$

$I_U^i$ is sometimes called the *conductance* of the set $U$ for the Markov chain associated to $\mathcal{L}^i$. For $t \in (0,1)$, the isoperimetric profile $\varphi_L^i(t)$ is defined by

$$\varphi_L^i(t) := \inf_{U \subset \xi_L \,:\, W^i(U) \leq tW^i(\xi_L)} I_U^i(\xi). \tag{1.11}$$

Cheeger's constant is defined by $\Phi_L^i := \varphi_L^i(\frac{1}{2})$.

THEOREM 1.1. *Assume Assumptions* (A1), (A2) *and* $\alpha < d$, $d \geq 2$. *For every* $i \in \{1, 2, 3\}$ *and every* $\varepsilon > 0$, *there exists* $\delta > 0$ *such that*, $\mathcal{P}$-*a.s.*,

$$\varphi_L^i(t) \geq \delta \min\left\{\frac{1}{L^\varepsilon}, \frac{1}{t^{1/d}L}\right\}, \qquad 0 < t \leq \frac{1}{2} \tag{1.12}$$

*for all $L$ sufficiently large. In particular, Cheeger's constant satisfies,* $\mathcal{P}$-*a.s.*,

$$\Phi_L^i \geq \delta L^{-1} \tag{1.13}$$

*for all $L$ sufficiently large.*

Recall that when $\xi$ coincides deterministically with $\mathbb{Z}^d$, the standard isoperimetric inequality implies that $\varphi_L^i(t) \geq \delta t^{-1/d} L^{-1}$ for any $0 < t \leq \frac{1}{2}$ (see, e.g., [25]). Theorem 1.1 shows that, up to scales that are smaller than any power of $L$, these estimates continue to hold in our three models.

The Poincaré constant $\gamma^i(L)$, $i = 1, 2, 3$, is defined by

$$\gamma^i(L) := \sup_{f:\xi_L \to \mathbb{R}} \frac{\text{Var}_i(f)}{\mathcal{E}_i(f)}, \tag{1.14}$$

where $\text{Var}_i(f)$ denotes the variance $\nu^i(f^2) - \nu^i(f)^2$ and $\mathcal{E}_i$ denotes the Dirichlet form

$$(1.15) \qquad \mathcal{E}_i(f) = \tfrac{1}{2} \sum_{x,y \in \xi_L} \nu^i(x) \mathcal{L}^i(x,y)(f(x) - f(y))^2.$$

$\gamma^i(L)$ is also known as the *relaxation time* and it coincides with the inverse of the spectral gap of the generator $\mathcal{L}^i$. Indeed, $-\mathcal{L}^i$ is a nonnegative definite matrix and its smallest nonzero eigenvalue $\lambda_1^i$ (the spectral gap) is given by $\lambda_1^i = 1/\gamma^i(L)$. The estimate (1.13) of Theorem 1.1 and an application of standard estimates (see Section 3) together imply that, assuming Assumptions (A1), (A2) and $\alpha < d$, there exists $C < \infty$ such that, $\mathcal{P}$-a.s.,

$$(1.16) \qquad \gamma^i(L) \leq CL^2$$

for all $L$ sufficiently large, for all $i = 1, 2, 3$. The Poincaré inequality (1.16) gives us some information on the speed with which the law of the random walk $X_t^i$ converges to the invariant distribution $\nu^i$. Namely, let

$$(1.17) \qquad H_t^i(x,y) = \sum_{n=0}^{\infty} \frac{t^n}{n!} (\mathcal{L}^i)^n(x,y)$$

denote the kernel of the random walk, that is, the probability that $X_t^i = y$ conditioned on $X_0^i = x$. We define the (uniform) *mixing time* $\tau^i(L)$ by

$$(1.18) \qquad \tau^i(L) = \inf\left\{ t > 0 \colon \sup_{x,y \in \xi_L} \frac{|H_t^i(x,y) - \nu^i(y)|}{\nu^i(y)} \leq \frac{1}{e} \right\}.$$

Well-known bounds (see, e.g., [25]) then allow us to estimate

$$(1.19) \qquad \tau^i(L) \leq \gamma^i(L)\left(1 + \log \frac{1}{\nu_*^i}\right),$$

where $\nu_*^i := \min_{x \in \xi_L} \nu^i(x)$. We shall see (Section 3.1) that under Assumptions (A1) and (A2), for $\alpha < d$ and $\varepsilon > 0$, we can deduce the uniform almost sure bound $\nu_*^i \geq \delta L^{-d-\varepsilon}$ for some constant $\delta > 0$. From (1.16), we then obtain that, for some $C < \infty$, $\mathcal{P}$-a.s.,

$$(1.20) \qquad \tau^i(L) \leq CL^2 \log L$$

for $L$ sufficiently large, $i = 1, 2, 3$.

Note that this estimate is only based on the bound (1.13) for Cheeger's constant. Theorem 1.1, however, shows that small sets can have a larger conductance. As first observed by Lovasz and Kannan [16], this fact can be used to obtain better bounds on the mixing time than (1.19). We will use refinements of this idea, recently obtained in [11, 21], to show that Theorem 1.1 implies the following improvement on (1.20).

THEOREM 1.2. *Assume Assumptions* (A1), (A2) *and* $\alpha < d$, $d \geq 2$. *There then exists* $C < \infty$ *such that, for every* $i = 1, 2, 3$, $\mathcal{P}$-*a.s.*,

$$\tau^i(L) \leq CL^2 \tag{1.21}$$

*for all L sufficiently large.*

Note that (1.21) is a strengthening of the estimate (1.16) since we always have (see, e.g., [25])

$$\gamma^i(L) \leq C\tau^i(L). \tag{1.22}$$

REMARK 1. Under suitable assumptions, one can show that the bound in Theorem 1.2 is tight up to constant factors. For instance, consider model 1 (i.e., $i = 1$) and assume that (1.5) holds $\mathcal{P}$-a.s. for $L$ sufficiently large. Set $f(x) = |x|$ in (1.14) and define $A_1 := \Lambda_{L/2}$ and $A_2 := \Lambda_{3L/4}$ so that

$$\begin{aligned}
\sum_{x,y \in \xi_L} (f(x) - f(y))^2 &\geq 2 \sum_{x \in \xi_L \cap A_1} \sum_{y \in \xi_L \cap A_2^c} (|x| - |y|)^2 \\
&\geq \xi(A_1)\xi(\Lambda_L \setminus A_2)\frac{L^2}{8}.
\end{aligned} \tag{1.23}$$

Suppose that there exists $c > 0$ such that, $\mathcal{P}$-a.s.,

$$\xi(A_1)\xi(\Lambda_L \setminus A_2) \geq cL^{2d} \tag{1.24}$$

for $L$ sufficiently large. Note that (1.24) always holds under Assumption (A1). Equations (1.23) and (1.24) then imply that $\mathrm{Var}_1(f) \geq cL^2$ for some constant $c > 0$. Suppose, further, that there exists $C < \infty$ such that, $\mathcal{P}$-a.s.,

$$\sum_{x,y \in \xi_L} e^{-|x-y|^\alpha}(|x| - |y|)^2 \leq CL^d \tag{1.25}$$

for $L$ sufficiently large. Then, $\mathcal{E}_1(f) \leq C$. It follows that there exists $\delta > 0$ such that, $\mathcal{P}$-a.s.,

$$\gamma^1(L) \geq \delta L^2 \tag{1.26}$$

for all $L$ sufficiently large. From (1.22), we derive the same estimate for $\tau^1(L)$. It is easy to check that (1.25) holds if, for example, Assumption (A2) is known to hold for all $\alpha > 0$. Indeed,

$$e^{-|x-y|^\alpha}(|x| - |y|)^2 \leq Ce^{-|x-y|^{\alpha/2}}.$$

The estimates in (1.16), Theorem 1.2 and (1.26) can be interpreted as the validity of diffusive behavior for the mixing time in the case $\alpha < d$, $d \geq 2$. Let us now turn to a discussion of the case $\alpha \geq d \geq 2$. We restrict our attention to homogeneous Poisson point processes $\mathcal{P} = \mathcal{P}_{*,\rho}$, but the same proof allows us to establish analogous results in a larger class of models including, for example, diluted lattices. Points (1) and (2) below give subdiffusive estimates on the Poincaré constant $\gamma^i(L)$ which, by (1.22), can be turned into estimates on $\tau^i(L)$. From points (2) and (3) below, we see that if $\alpha = d$, for $\mathcal{P} = \mathcal{P}_{*,\rho}$, there is a transition from subdiffusive to diffusive behavior as the intensity $\rho$ is increased.

THEOREM 1.3. *For $i = 1, 2, 3$, we have the following:*

(1) $\alpha > d$. *For every $\rho > 0$, for any $\delta > 0$,*

$$\gamma^i(L) \geq L^{1/\delta}, \tag{1.27}$$

$\mathcal{P}_{*,\rho}$-*a.s., for $L$ sufficiently large.*

(2) $\alpha = d$ *(small density). For every $\delta > 0$, there exists $\rho_0 > 0$ such that, for any $0 < \rho \leq \rho_0$,*

$$\gamma^i(L) \geq L^{1/\delta}, \tag{1.28}$$

$\mathcal{P}_{*,\rho}$-*a.s., for $L$ sufficiently large.*

(3) $\alpha = d$ *(high density). There exist $C < \infty$ and $\rho_1 < \infty$ such that, for any $\rho \geq \rho_1$,*

$$\tau^i(L) \leq CL^2, \tag{1.29}$$

$\mathcal{P}_{*,\rho}$-*a.s., for $L$ sufficiently large.*

The rest of the paper is organized as follows. In Section 2, we prove Theorem 1.1. Theorem 1.2 is proved in Section 3. In Section 4, we prove Theorem 1.3. In Section 5, we discuss examples of point processes satisfying the assumptions required in the main statements. In Appendix A, we prove that Assumption (A2) is satisfied by the homogeneous Poisson point process, while in Appendix B, we collect some preliminary results on Bernoulli site percolation.

**2. Isoperimetric inequalities.**

2.1. *The basic construction.* We shall partition the space $\mathbb{R}^d$ according to three different scales: $K$, $L^\varepsilon$, $C_W (\log L)^{1/d}$, where $\varepsilon > 0$ is a small constant and $K, C_W$ are suitably large constants (independent of $L$). Accordingly, we will denote by $B_x$ the cubes of side $K$ (as in the previous section), by $V_x$ the cubes of side $L^\varepsilon$ and by $W_x$ the cubes of side $C_W (\log L)^{1/d}$. These are all assumed to be of the form $ax + [0, a)^d$, $x \in \mathbb{Z}^d$, for $a = K, L^\varepsilon, C_W (\log L)^{1/d}$, respectively. At the cost of replacing $a$ with $(L/2)/\lfloor (L/2)/a \rfloor = a(1 + o(1))$, with some abuse of notation, we assume that the box $\Lambda_L$ is partitioned by the $a$-cubes.

2.1.1. *The cluster of grey cubes.* Let us look at the partition into $K$-cubes $B_x$. As in the Introduction, we call a cube $B_x$ such that $\xi(B_x) \geq 1$ *good*. On the same probability space of the point process, we consider the independent Bernoulli process which assigns to a box $B_x$ the color *grey* with probability $p$ and the color *white* with probability $1 - p$. From our Assumption (A1), we may assume that whenever a cube $B_x$ is grey, it is also good. A collection of cubes is said to be *connected* if any two cubes $B_x, B_y$ belonging to it can be joined by a path of adjacent cubes, where two cubes $B_z$ and $B_{z'}$ are said to be *adjacent* if their centers are at distance $K$. We denoted by $\mathcal{C}_L$ the largest connected component of grey cubes $B_x$ such that $B_x \subset \Lambda_L$; see Figure 1. This is well defined since, with probability 1, there is eventually a unique cluster with maximal cardinality if $p < 1$ is sufficiently large, as will be assumed below (see, e.g., [13] or Appendix B). We refer to $\mathcal{C}_L$ as the *cluster of grey cubes*.

2.1.2. *Density of the cluster of grey cubes.* We shall need the following consequence of Assumption (A1). Let $\{V_x\}_{x \in J_L}$ denote the partition of $\Lambda_L$ by means of $L^\varepsilon$-cubes. There then exists $\delta > 0$ such that if $K$ is sufficiently large, then, $\mathcal{P}$-a.s.,

$$(2.1) \qquad |\mathcal{C}_L \cap V_x| \geq \delta L^{\varepsilon d} \qquad \forall x \in J_L$$

for all $L$ sufficiently large. Here, $|\mathcal{C}_L \cap V_x|$ stands for the volume of the intersection as a subset of $\mathbb{R}^d$.

Due to Assumption (A1), the grey $K$-cubes correspond to a Bernoulli site percolation with parameter $p$ that can be taken arbitrarily close to 1 by choosing $K$ large. Hence, it is enough to check the above result for Bernoulli site percolation with $p < 1$ sufficiently large. This is done in Appendix B.

Another consequence of Assumption (A1) is that every $C_W (\log L)^{1/d}$-cube $W_x$ in the partition of $\Lambda_L$ must contain at least one grey $K$-cube $B_z$ and hence at least one point of $\xi_L$. Indeed, the probability that there exists one cube $W_x \subset \Lambda_L$ containing no grey $K$-cube is bounded from above by

$$(2.2) \qquad \left(\frac{L}{C_W \log L}\right)^d (1-p)^{(C_W^d/K^d)\log L},$$

which is summable in $L \in \mathbb{N}$ if $C_W$ is sufficiently large. The Borel–Cantelli lemma then shows that, almost surely, every $W_x \subset \Lambda_L$ eventually contains at least one grey $K$-cube.

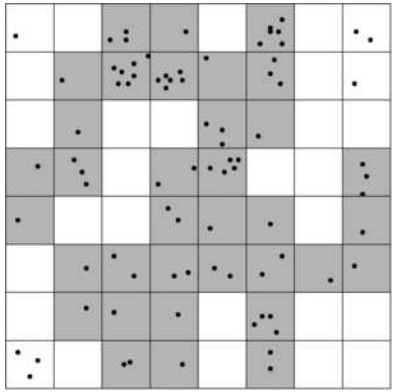

FIG. 1. *A realization $\xi_L$ of the point process in $\Lambda_L$ and the cluster of grey cubes $\mathcal{C}_L$.*

2.1.3. *Isoperimetric inequality for $\mathcal{C}_L$.* As shown in [2] and [20], if $p$ is close to 1, it is not very hard to establish good isoperimetric bounds for the percolation cluster. We shall now state these estimates explicitly in our context.

For any collection $A \subset \mathcal{C}_L$ of $K$-cubes, we define $\partial A$ as the collection of $K$-cubes $B_y$ such that $B_y \in \mathcal{C}_L \setminus A$ and $B_y$ is adjacent to some $K$-cube belonging to $A$. $|A|$ and $|\partial A|$ will denote their respective volumes as subsets of $\mathbb{R}^d$.

LEMMA 2.1. *For suitably large values of the constant $K$, the following holds. There exist positive constants $\kappa, C$ such that, $\mathcal{P}$-almost surely, for $L$ sufficiently large, every collection $A \subset \mathcal{C}_L$ such that $\frac{1}{2}|\mathcal{C}_L| \geq |A| \geq C(\log L)^{d/(d-1)}$ satisfies*

$$(2.3) \qquad \frac{|\partial A|}{|A|} \geq \kappa \frac{1}{|A|^{1/d}}.$$

Since grey cubes form a Bernoulli process with parameter $p$ that can be taken to be close to 1 (by choosing $K$ suitably large), the proof of Lemma 2.1 follows from Lemma 2.6 of [2].

When the collection $A$ is such that $0 < |A| \leq C(\log L)^{d/(d-1)}$, we simply observe that since $\partial A$ is nonempty, $|\partial A| \geq K^d$ and therefore

$$(2.4) \qquad \frac{|\partial A|}{|A|} \geq \frac{K^d}{C(\log L)^{d/(d-1)}}.$$

Together with Lemma 2.1, this shows that for any $A \subset \mathcal{C}_L$ such that $\frac{1}{2}|\mathcal{C}_L| \geq |A| > 0$, we have

$$(2.5) \qquad \frac{|\partial A|}{|A|} \geq \kappa \min\left\{\frac{1}{|A|^{1/d}}, \frac{1}{(\log L)^{d/(d-1)}}\right\}$$

for a suitable constant $\kappa > 0$.

2.2. *Proof of Theorem* 1.1. We are given an arbitrary set $U \subset \xi_L$ with $W^i(U) \leq \frac{1}{2} W^i(\xi_L)$ and we have to estimate

$$I_U^i = \frac{1}{W^i(U)} \sum_{x \in U, y \in U^c} e^{-|x-y|^\alpha}, \tag{2.6}$$

where $U^c$ stands for the complement w.r.t. $\xi_L$, that is, $U^c = \xi_L \setminus U$.

We denote by $S = S(U)$ the set of cubes $B_x \in \mathcal{C}_L$ which intersect $U$ (see Figure 2):

$$S = \{B_x \in \mathcal{C}_L : B_x \cap U \neq \varnothing\}. \tag{2.7}$$

Also, we let $T = T(U)$ denote the set of cubes $B_x \in \mathcal{C}_L$ which intersect $U^c$:

$$T = \{B_x \in \mathcal{C}_L : B_x \cap U^c \neq \varnothing\}. \tag{2.8}$$

Roughly speaking, the idea of the proof will be to exploit as much as possible the isoperimetric estimates for the region $S$ as a subset of $\mathcal{C}_L$ (Lemma 2.1). We shall consider two separate cases, according to whether $W^i(U) > \gamma |S|$ or $W^i(U) \leq \gamma |S|$, where $\gamma$ is a suitably large constant to be fixed below. If $W^i(U) > \gamma |S|$, we show that, in this case, the set $U$ contains islands, whose number is [neglecting $O(L^\varepsilon)$ corrections] proportional to $W^i(U)$ and each of which contributes at least $L^{-\varepsilon}$ to the numerator in $I_U^i$. If $W^i(U) \leq \gamma |S|$ and the set $S$ is not too large, we can essentially rely on Lemma 2.1. Finally, it will remain to discuss the case when $S$ is large in $\mathcal{C}_L$ (e.g., $S = \mathcal{C}_L$). In this case, either the set $T$ is large as well, in which case we must have $S \cap T$ large and this produces a large numerator in $I_U$, or $T$ is small compared to $W^i(U^c)$, in which case things will be handled as in the case $W^i(U) > \gamma |S|$ discussed above, by switching from $U$ to $U^c$.

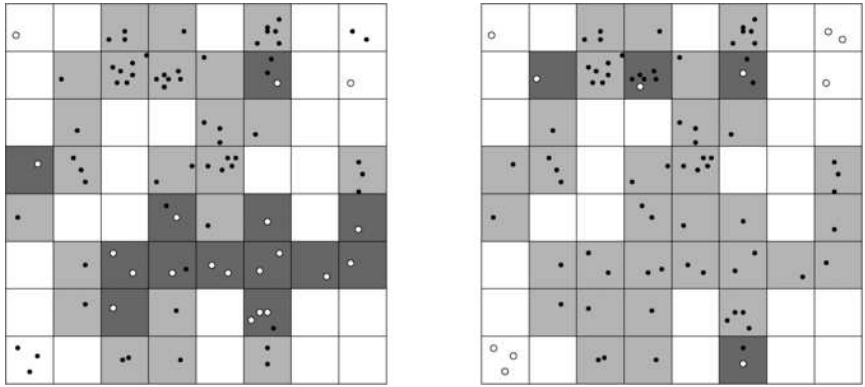

FIG. 2. *Consider the same realizations $\xi_L, \mathcal{C}_L$ as in Figure 1. To mark the chosen set $U$, we use $\circ$ for points in $U$ and $\bullet$ for points in $U^c$. Two different choices of the set $U$ are depicted above. In each case, the set $S \subset \mathcal{C}_L$ is given by the union of the darker cubes.*

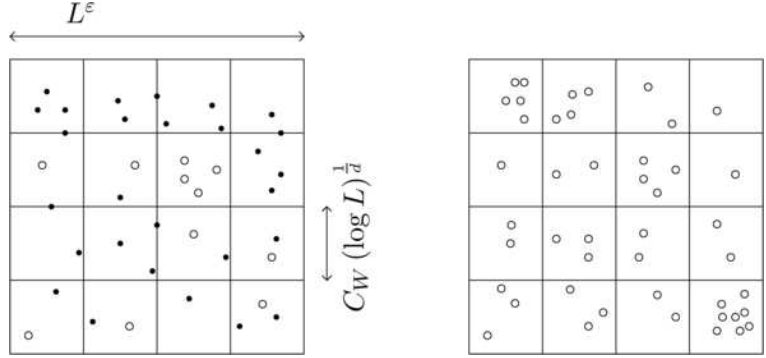

Fig. 3. *Two special cubes, strong (left) and weak (right).*

2.2.1. *The case $W^i(U) > \gamma|S|$.*

LEMMA 2.2. *For every $\varepsilon > 0$, there exists a constant $\gamma < \infty$ such that, $\mathcal{P}$-a.s., for $L$ sufficiently large, we have*

$$I_U^i \geq L^{-\varepsilon} \tag{2.9}$$

*for all $U \subset \xi_L$ such that $W^i(U) > \gamma|S|$.*

PROOF. Observe that, from definitions (1.3) and (1.9), we have $R_A \geq W^i(A)$ for every bounded set $A$, for every $i = 1, 2, 3$. Recall the definition of the $L^\varepsilon$-cubes $V_x$. Let $m$ denote the number of cubes $V_x$ such that $U \cap V_x \neq \varnothing$. Such $L^\varepsilon$-cubes are called *special*. Since $R_{V_x}(\xi) \geq W^i(V_x)$, by Assumption (A2), we may assume that $W^i(V_x) \leq CL^{\varepsilon d}$ for all $V_x$ in the partition of $\Lambda_L$, so

$$m \geq \frac{W^i(U)}{CL^{\varepsilon d}} \tag{2.10}$$

for some constant $C < \infty$. Special cubes $V_x$ are classified as either *strong* and *weak* according to whether, besides points of $U$, they also contain points of $U^c$. Namely, a special cube $V$ is called strong if $V \cap U^c \neq \varnothing$ and weak if $V \cap U^c = \varnothing$; see Figure 3.

Suppose $V_w$ is a given weak special cube. Since $V_w \cap U^c = \varnothing$, we must have

$$S \cap V_w = \mathcal{C}_L \cap V_w, \tag{2.11}$$

where, with some abuse of notation we identify a collection of cubes with their union. Indeed, if $B_y$ is a $K$-cube with $B_y \in V_w \cap \mathcal{C}_L$, then we must have $B_y \cap \xi_L \neq \varnothing$ and $B_y \cap U^c = \varnothing$, hence $B_y \cap U \neq \varnothing$ and therefore $B_y \in S$. Thus, by the density property (2.1), we have

$$|S \cap V_w| \geq \delta L^{\varepsilon d} \tag{2.12}$$

for some constant $\delta > 0$. Letting $m_1$ denote the number of weak special cubes in $\Lambda_L$ and summing the contributions (2.12) of each weak special cube $V_w$, we obtain

$$|S| \geq \delta m_1 L^{\varepsilon d}. \tag{2.13}$$

We now observe that there are more strong special cubes than weak special cubes, that is, that $m_1 < m/2$. Indeed, otherwise, from (2.10) and (2.13), we would obtain

$$|S| \geq \delta \frac{m}{2} L^{\varepsilon d} \geq \frac{\delta}{2C} W^i(U), \tag{2.14}$$

which contradicts $W^i(U) > \gamma |S|$ if $\gamma$ is sufficiently large.

Let $V_s$ be a given strong special cube. We claim that

$$\sum_{\substack{x \in U \cap V_s, \\ y \in U^c \cap V_s}} e^{-|x-y|^\alpha} \geq e^{-c(\log L)^{\alpha/d}} \tag{2.15}$$

for $c := (2\sqrt{d} C_W)^\alpha$. To establish (2.15), it is sufficient to show that we can find two vertices $x, y \in V_s$ such that $x \in U$, $y \in U^c$ and $|x-y| \leq 2\sqrt{d} C_W (\log L)^{1/d}$.

The two points $x$ and $y$ with the above property can be found as follows. Let $\{W_i\}$ denote the collection of $C_W (\log L)^{1/d}$-cubes $W_i$ such that $W_i \subset V_s$. Since $V_s \subset \Lambda_L$, from (2.2), we know that each of the cubes $W_i$ contains at least one point of $\xi_L$. Moreover, by the definition of a strong special cube, we know that $V_s$ contains a point $z \in U$ and a point $z' \in U^c$. Let $z = z_0, z_1, \ldots, z_n = z'$ denote a path joining the two points $z, z'$, such that:

(1) $z_i \in \xi_L$ for every $i = 1, \ldots, n-1$;
(2) $|z_i - z_{i+1}| \leq 2\sqrt{d} C_W (\log L)^{1/d}$ for every $i = 1, \ldots, n-1$.

Such a path exists, since each of the cubes $W_i$ contains at least one point of $\xi_L$. Let $k = \min\{i \geq 1 : z_i \in U^c\}$. The needed points $x, y$ are obtained by setting $x := z_{k-1}$, $y := z_k$. This completes the proof of (2.15).

From (2.14), we know that there are at least $m/2$ strong special cubes. Restricting to strong special cubes $V_s$ and using (2.15), we can then estimate

$$\sum_{x \in U, y \in U^c} e^{-|x-y|^\alpha} \geq \sum_{V_s} \sum_{\substack{x \in U \cap V_s, \\ y \in U^c \cap V_s}} e^{-|x-y|^\alpha} \geq \frac{m}{2} e^{-c(\log L)^{\alpha/d}}. \tag{2.16}$$

Since $\alpha < d$, we can estimate $e^{-c(\log L)^{\alpha/d}} \geq L^{-\varepsilon}$. Using (2.16) and (2.10), we see that for every $\varepsilon > 0$, for any $L$ sufficiently large,

$$I_U^i = \frac{1}{W^i(U)} \sum_{x \in U, y \in U^c} e^{-|x-y|^\alpha} \geq \frac{e^{-c(\log L)^{\alpha/d}}}{2C L^{\varepsilon d}} \geq L^{-2\varepsilon d}. \tag{2.17}$$

This completes the proof of Lemma 2.2. □

REMARK 2. We point out that the proof of Lemma 2.2 never used the fact that $W^i(U) \leq \frac{1}{2}W^i(\xi_L)$. This will be needed in the remainder of the argument below.

REMARK 3. Also, we note that (2.17) is the only piece of the proof using the assumption $\alpha < d$. In fact, the same estimate as in (2.17) would hold in the case $\alpha = d$ if one could choose the constant $c$ arbitrarily small. Since $c = (2\sqrt{d}C_W)^\alpha$, this can be achieved by taking $C_W$ small. However, the constant $C_W$ must be sufficiently large in order to guarantee that, almost surely, all $C_W(\log L)^{1/d}$-cubes $W_x \subset \Lambda_L$ are occupied by at least one point of $\xi$, for $L$ sufficiently large; see (2.2). These observations will be used in the proof of Theorem 1.3.

2.2.2. *The case $W^i(U) \leq \gamma|S|$.*

LEMMA 2.3. *Let $\varepsilon, \gamma$ be the constants appearing in Lemma 2.2. There then exists $\delta > 0$ such that, $\mathcal{P}$-a.s., for all $L$ sufficiently large, we have*

(2.18) $$I_U^i \geq \delta \min\{L^{-\varepsilon}, W^i(U)^{-1/d}\}$$

*for all $U \subset \xi_L$ such that $W^i(U) \leq \gamma|S|$ and $W^i(U) \leq \frac{1}{2}W^i(\xi_L)$.*

PROOF. From our basic construction, any $K$-cube in $S$ has at least one point in $U$ and any $K$-cube in $\mathcal{C}_L \setminus S$ has at least one point in $U^c$. Therefore,

(2.19) $$\sum_{x \in U, y \in U^c} e^{-|x-y|^\alpha} \geq \delta \max\{|\partial S|, |\partial(\mathcal{C}_L \setminus S)|\},$$

where $\delta = \delta(d, K, \alpha)$ is a positive constant.

We now prove our claim under the assumption that $|S| \leq a|\mathcal{C}_L|$ for some given constant $a \in (\frac{1}{2}, 1)$. We will remove this restriction afterward. In this case, we have

$$W^i(U) \leq \gamma|S| \leq \gamma \min\left\{|S|, \frac{a}{1-a}|\mathcal{C}_L \setminus S|\right\}.$$

Setting $S^\# = S$ if $|S| \leq \frac{1}{2}|\mathcal{C}_L|$ and $S^\# = \mathcal{C}_L \setminus S$ if $|S| > \frac{1}{2}|\mathcal{C}_L|$, from (2.19), we obtain

(2.20) $$I_U^i \geq \frac{\delta(1-a)}{\gamma a} \frac{|\partial S^\#|}{|S^\#|}.$$

Using (2.5), we conclude that

(2.21) $$I_U^i \geq \frac{\kappa\delta(1-a)}{\gamma a} \min\left\{\frac{1}{|S^\#|^{1/d}}, \frac{1}{(\log L)^{d/(d-1)}}\right\}.$$

Next, we claim that

$$|S| \leq CW^i(U). \tag{2.22}$$

For $i = 1$, we have $W^1(U) = \#(U)$ and, since there is at least one point of $U$ in each cube of $S$, we have $|S| \leq K^d \#(U)$. For $i = 2$, observe that for every $x \in U \cap S$, there exists $y \in \xi_L \setminus \{x\}$ such that $|x - y| \leq C$ for a constant $C = C(K, d)$. Indeed, let $B$ denote the $K$-cube such that $x \in B$ and let $B'$ denote a cube adjacent to $B$ in $\mathcal{C}_L$. Then, by construction, $B'$ contains at least one point $y \in \xi_L$. Therefore, $w_x^2 \geq e^{-|x-y|^\alpha} \geq e^{-C} =: c$ for every $x \in U \cap S$. It follows that $W^2(U) \geq W^2(S \cap U) \geq cK^{-d}|S|$. Finally, for $i = 3$, simply use the fact that $W^3(U) \geq W^1(U)$. This proves (2.22).

Since $|S^\#| \leq |S|$, using (2.22) and (2.21), we arrive at the bound

$$I_U^i \geq \frac{\kappa\delta(1-a)}{\gamma a} \min\left\{\frac{1}{CW^i(U)^{1/d}}, \frac{1}{(\log L)^{d/(d-1)}}\right\}. \tag{2.23}$$

This proves the claim (2.18) under the assumption $|S| \leq a|\mathcal{C}_L|$.

We must now remove the latter restriction. In particular, nothing prevents our set $S$ from coinciding with $\mathcal{C}_L$. Suppose, then, that $|S| > a|\mathcal{C}_L|$. Let $T$ denote the set defined by (2.8) and assume that $|T| \geq 2(1-a)|\mathcal{C}_L|$. Then, $|\mathcal{C}_L \setminus T| \leq (2a-1)|\mathcal{C}_L|$ and $|S| \leq |S \cap T| + (2a-1)|\mathcal{C}_L|$. Therefore, in this case,

$$|S \cap T| \geq (1-a)|\mathcal{C}_L|. \tag{2.24}$$

Clearly, any cube $B_y \in S \cap T$ contains at least one point of $U$ and at least one point of $U^c$. Therefore,

$$\sum_{x \in U, y \in U^c} e^{-|x-y|^\alpha} \geq \delta_1 |S \cap T|, \tag{2.25}$$

where $\delta_1 = \delta_1(d, K, \alpha)$ is a positive constant. Now, $W^i(U) \leq \gamma|S| \leq \gamma|\mathcal{C}_L|$ and, therefore, (2.24) and (2.25) imply that

$$I_U^i \geq \frac{\delta_1}{\gamma}(1-a). \tag{2.26}$$

The bound (2.26) was obtained using $|T| \geq 2(1-a)|\mathcal{C}_L|$. If, on the contrary, $|T| < 2(1-a)|\mathcal{C}_L|$, then observe that, as in (2.22), we have the bound $W^i(\xi_L) \geq c|\mathcal{C}_L|$ for some positive constant $c$ and therefore

$$W^i(U^c) \geq \frac{1}{2}W^i(\xi_L) \geq \frac{c}{2}|\mathcal{C}_L| > \frac{c}{4}\frac{1}{(1-a)}|T|. \tag{2.27}$$

Therefore, $W^i(U^c) > \gamma|T|$ if $a$ is sufficiently close to 1, where $\gamma$ was fixed in Lemma 2.2.

Now, since $W^i(U) \leq W^i(U^c)$ [because $W^i(U) \leq \frac{1}{2}W^i(\xi_L)$, by assumption], we have $I_U^i \geq I_{U^c}^i$. Moreover, the collection of $K$-cubes $T$ is, for the set $U^c$, exactly what $S$ is for $U$; see (2.7) and (2.8). As discussed in Remark 2, we can repeat the argument of Lemma 2.2 with $U$ replaced by $U^c$ and $S$ replaced by $T$ since that argument applies, despite the fact that we now have $W^i(U^c) \geq \frac{1}{2}W^i(\xi_L)$. Using the bound $W^i(U^c) > \gamma|T|$, we can therefore estimate

$$I_U^i \geq I_{U^c}^i \geq L^{-\varepsilon}. \tag{2.28}$$

This completes the proof of Lemma 2.3. □

2.2.3. *Conclusion.* To finish the proof of Theorem 1.1, we need only gather the estimates in Lemmas 2.2 and 2.3 to obtain that, for every $i = 1, 2, 3$, for arbitrary $\varepsilon > 0$, there exists $c > 0$ such that, $\mathcal{P}$-a.s., for all $L$ sufficiently large, we have

$$I_U^i \geq c \min\{L^{-\varepsilon}, W^i(U)^{-1/d}\} \tag{2.29}$$

for any set $U \subset \xi_L$ with $W^i(U) \leq \frac{1}{2}W^i(\xi_L)$. Passing to the function $\varphi_L^i$ defined in (1.11), we arrive at the almost sure estimate

$$\varphi_L^i(t) \geq c \min\left\{\frac{1}{L^{\varepsilon}}, \frac{1}{t^{1/d}W^i(\xi_L)^{1/d}}\right\}. \tag{2.30}$$

Since, from Assumption (A2), we know that $W^i(\xi_L) \leq cL^d$ almost surely, (2.30) immediately implies the estimate in Theorem 1.1.

REMARK 4. For further applications related to the Palm distribution (see Section 5.5), we point out that the proof given above can easily be adapted to show that for every $\varepsilon > 0$, there exists $\delta > 0$ such that, $\mathcal{P}$-a.s., for all $0 < t \leq 1/2$,

$$\inf_{z \in \Lambda_1} \inf_{U \subset \xi_{z,L-1} : W^i(U) \leq tW^i(\Lambda_{z,L-1})} \frac{1}{W^i(U)} \sum_{x \in U, y \in \xi_{z,L-1} \setminus U} e^{-|x-y|^{\alpha}}$$

$$\geq \delta \min\left\{\frac{1}{L^{\varepsilon}}, \frac{1}{t^{1/d}L}\right\} \tag{2.31}$$

for $L$ sufficiently large, where

$$\xi_{z,L-1} := \xi \cap \Lambda_{z,L-1}, \qquad \Lambda_{z,L-1} := \Lambda_{L-1} + z.$$

We sketch below the argument needed to derive (2.31).

Fix $z \in \Lambda_1$, consider $U \subset \xi_{z,L-1}$ with $W^i(U) \leq \frac{1}{2}W^i(\xi_{z,L-1})$ and set $U^c := \xi_{z,L-1} \setminus U$. Replace, in the previous proof, $\mathcal{C}_L$ with $\mathcal{C}_{L-2}$ and keep all the remaining notation (up to the above substitutions). Note that $\xi_{z,L-1} \subset \Lambda_{z,L-1} \subset \Lambda_L$ and that, due to (1.5),

$$c_1 L^d \leq \#(\xi_{L-2}) \leq \#(\xi_{z,L-1}) \leq \#(\xi_L) \leq c_2 L^d \tag{2.32}$$

for suitable positive constants $c_1, c_2$ which do not depend on $z$. Moreover, observe that (2.1) still holds with $\mathcal{C}_L$ replaced with $\mathcal{C}_{L-2}$ and that, by choosing $C_W$ sufficiently large, every $C_W(\log L)^{1/d}$-cube included in $\Lambda_{L-2}$ must contain a point of $\xi$ and hence a point of $\xi_{z,L-1}$ [see the derivation of (2.2)].

Consider the case $W^i(U) > \gamma|S|$ of Section 2.2.1. All of the arguments before (2.15) remain valid, but the explanation of (2.11) is now as follows: every $K$-cube in $\mathcal{C}_{L-2}$ contains a point in $\xi_{L-2} \subset \xi_{z,L-1}$ and, in particular, if $B_y$ is a $K$-cube with $B_y \in V_w \cap \mathcal{C}_{L-2}$, it must be $B_y \cap \xi_{z,L-1} \neq \varnothing$ and $B_y \cap U^c = \varnothing$, hence $B_y \cap U \neq \varnothing$ and therefore $B_y \in S$. In order to establish (2.15), we need to show that each strong special cube $V_s$ contains points $x, y$ such that $x \in U$, $y \in U^c$ and $|x - y| \leq c(d)C_W(\log L)^{1/d}$. The arguments given in the proof of Theorem 1.1 need to be slightly modified by defining $\{W_i\}$ as the family of $C_W(\log L)^{1/d}$-cubes included in $\Lambda_{L-2}$ and observing that each $W_i$ must contain a point of $\xi_{z,L-1}$, while the set $V_s \setminus \Lambda_{L-2}$ is very thin since it is included in $\Lambda_L \setminus \Lambda_{L-2}$. Having (2.15), (2.17) is derived as in the proof of Theorem 1.1. Neither the case $W^i(U) \leq \gamma|S|$ of Section 2.2.1 nor the conclusions of Section 2.2.3 need any additional modification.

**3. Upper bounds on mixing times.** We are going to prove Theorem 1.2. We start by recalling the so-called *spectral profile function* and its use in bounding mixing times [11]. For models $i = 1, 2, 3$, for any $U \subset \xi_L$, define

$$\lambda^i(U) = \inf_{f \in c_0^+(U)} \frac{\mathcal{E}_i(f)}{\text{Var}_i(f)}, \tag{3.1}$$

where $c_0^+(U)$ denotes the set of functions $f : \xi_L \to \mathbb{R}$ such that $f \geq 0$ and $f(x) = 0$ for all $x \in U^c = \xi_L \setminus U$. Recall that $\nu_*^i := \min_{x \in \xi_L} \nu^i(x)$. The spectral profile function $\Lambda^i : [\nu_*^i, \infty) \to \mathbb{R}$ is defined by

$$\Lambda^i(r) = \inf_{U \subset \xi_L \,:\, \nu_*^i \leq \nu^i(U) \leq r} \lambda^i(U). \tag{3.2}$$

The main result in [11] can be restated as follows.

LEMMA 3.1. *For every $i = 1, 2, 3$, the mixing time $\tau^i(L)$ satisfies*

$$\tau^i(L) \leq 2 \int_{4\nu_*^i}^{4e} \frac{dr}{r\Lambda^i(r)}. \tag{3.3}$$

Lemma 3.1 is contained in Theorem 1.1 of [11], which is a general result for continuous-time Markov chains with finite state space. The next step is to bound the spectral profile in terms of our isoperimetric profile $\varphi_L^i$. Namely, we need a bound of the form

$$\Lambda^i(r) \geq \delta \varphi_L^i(r)^2, \qquad r \in [\nu_*^i, 1). \tag{3.4}$$

We shall obtain this bound for models 2 and 3. There is a technical difficulty in obtaining the same estimate for model 1, so we shall treat this case separately.

3.1. *Models* 2 *and* 3. The estimate (3.4) is derived in Lemma 2.4 of [11] for the case of Markov chains with generator $\mathcal{L}$ of the form $K-1$, where $K$ is a stochastic matrix. Note that $\mathcal{L}^i$ is of this form for models $i = 2, 3$. We shall give the details in the next lemma for the reader's convenience.

LEMMA 3.2. *For $i = 2, 3$, we have*

$$\Lambda^i(r) \geq \tfrac{1}{2} \varphi_L^i(r)^2, \qquad r \in [\nu_*^i, 1). \tag{3.5}$$

*Moreover,*

$$\gamma^i(L) \leq 2(\Phi_L^i)^{-2}, \qquad i = 2, 3. \tag{3.6}$$

PROOF. Let $f : \xi_L \to \mathbb{R}_+$ be a nonnegative function. Set $F_t = \{x \in \xi_L : f(x) \geq t\}$. Then,

$$\nu^i(f) = \sum_{x \in \xi_L} \nu^i(x) f(x) = \int_0^\infty \nu^i(F_t) \, dt. \tag{3.7}$$

Set $Q^i(x, y) := \nu^i(x) \mathcal{L}^i(x, y)$ and $Q^i(U, V) := \sum_{x \in U, y \in V} Q^i(x, y)$ for any $U, V \subset \xi_L$. Note that, for any $i = 1, 2, 3$, $Q^i(x, y) = Q^i(y, x)$. We then have

$$\begin{aligned}
\sum_{x, y \in \xi_L} &|f(x) - f(y)| Q^i(x, y) \\
&= 2 \sum_{x, y \in \xi_L : f(x) > f(y)} [f(x) - f(y)] Q^i(x, y) \\
&= 2 \sum_{x, y \in \xi_L : f(x) > f(y)} Q^i(x, y) \int_0^\infty 1_{\{f(x) \geq t > f(y)\}} \, dt \\
&= 2 \int_0^\infty Q^i(F_t, F_t^c) \, dt.
\end{aligned} \tag{3.8}$$

Now, recall that

$$I_{F_t}^i = \frac{Q^i(F_t, F_t^c)}{\nu^i(F_t)}.$$

If $f(x) = 0$ for all $x \in U^c$, then $F_t \subset U$ for all $t > 0$ and therefore $I_{F_t}^i \geq \varphi_L^i(\nu^i(U))$. From (3.8) and (3.7), we have thus obtained

$$\begin{aligned}
\sum_{x, y \in \xi_L} |f(x) - f(y)| Q^i(x, y) &\geq 2 \varphi_L^i(\nu^i(U)) \int_0^\infty \nu^i(F_t) \, dt \\
&= 2 \varphi_L^i(\nu^i(U)) \nu^i(f)
\end{aligned} \tag{3.9}$$

for all $f \in c_0^+(U)$. For any such $f$, we may apply (3.9) to $f^2 \in c_0^+(U)$ to obtain

$$
\begin{aligned}
2\varphi_L^i(\nu^i(U))\nu^i(f^2) &\leq \sum_{x,y \in \xi_L} |f(x)^2 - f(y)^2| Q^i(x,y) \\
&= \sum_{x,y \in \xi_L} |f(x) - f(y)||f(x) + f(y)| Q^i(x,y) \\
&\leq \left(\sum_{x,y \in \xi_L} (f(x) - f(y))^2 Q^i(x,y)\right)^{1/2} \\
&\quad \times \left(\sum_{x,y \in \xi_L} (f(x) + f(y))^2 Q^i(x,y) 1_{x \neq y}\right)^{1/2} \\
&\leq (2\mathcal{E}_i(f))^{1/2} \left(4 \sum_{x,y \in \xi_L} f(x)^2 Q^i(x,y) 1_{x \neq y}\right)^{1/2},
\end{aligned}
$$

where we use Schwarz' inequality and the symmetry of $Q^i$.

Now, observe that, for $i = 2, 3$, we have $\sum_{y \in \xi_L : y \neq x} \mathcal{L}^i(x,y) \leq 1$ for every $x$, so

$$\sum_{x,y \in \xi_L} f(x)^2 Q^i(x,y) 1_{x \neq y} \leq \nu^i(f^2)$$

in these cases. This shows that, for $i = 2, 3$,

$$(3.10) \qquad 2\varphi_L^i(\nu^i(U))\nu^i(f^2) \leq (2\mathcal{E}_i(f))^{1/2}(4\nu^i(f^2))^{1/2}.$$

Therefore, from (3.1) and (3.10), we see that

$$(3.11) \qquad \lambda^i(U) \geq \inf_{f \in c_0^+(U)} \frac{\mathcal{E}_i(f)}{\nu^i(f^2)} \geq \frac{1}{2}\varphi_L^i(\nu^i(U))^2.$$

Returning to the profile functions, we have $\Lambda^i(r) \geq \frac{1}{2}\varphi_L^i(r)^2$ for any $r \in [\nu_*^i, 1)$. In a similar way (see Remark 2.1 in [11]), one proves that the Poincaré constants $\gamma^i(L)$ satisfy (3.6) for $i = 2, 3$. $\square$

REMARK 5. The only difficulty in obtaining the same type of estimates for $i = 1$ is that the sum $\sum_{y \in \xi_L : y \neq x} \mathcal{L}^1(x,y)$ cannot be given a uniform upper bound. This is where the third model ($i = 3$) becomes useful; see Section 3.2 below.

We are now able to prove Theorem 1.2 for $i = 2, 3$. When $r > \frac{1}{2}$, we can use the simple bound $\Lambda^i(r) \geq (\gamma^i(L))^{-1} \geq \frac{1}{2}(\Phi_L^i)^2$ (the first bound follows from the definitions, while the latter is implied by Lemma 3.2). Therefore, as a corollary of Lemma 3.2 and Lemma 3.1, we obtain that for any $L \in \mathbb{N}$, the mixing times satisfy

$$\tau^i(L) \leq 4 \int_{4\nu_*^i}^{4e} \frac{dt}{t\bar{\varphi}_L^i(t)^2}, \tag{3.12}$$

where $\bar{\varphi}_L^i(t) := \varphi_L^i(t)$ for $t \leq \frac{1}{2}$ and $\bar{\varphi}_L^i(t) := \Phi_L^i$ for $t > \frac{1}{2}$.

To complete the proof, we take $\varepsilon = \frac{1}{2}$ in Theorem 1.1 and observe that, neglecting multiplicative constants, the integral in (3.12) can be bounded from above by

$$\int_{4\nu_*^i}^{L^{-d/2}} \frac{dt}{tL^{-1}} + \int_{L^{-d/2}}^{1/2} \frac{dt}{t^{1-2/d}L^{-2}} + \int_{1/2}^{4e} \frac{dt}{tL^{-2}}. \tag{3.13}$$

Note that $\nu_*^i \geq CL^{-d-\varepsilon}$. Indeed, from Assumption (A1), we know that every cube $W_x$ contains at least one point of $\xi_L$ [see (2.2)] and therefore, for every $x \in \xi_L$, we have

$$w_x^3 \geq w_x^2 \geq e^{-C(\log L)^{\alpha/d}} \geq L^{-\varepsilon}$$

for $\alpha < d$. On the other hand, from Assumption (A2), we know that $W^{2,3}(\xi_L) \leq CL^d$. It follows that $\nu^{2,3}(x) \geq CL^{-d-\varepsilon}$ for every $x \in \xi_L$. This shows that the first term in (3.13) contributes at most $O(L \log L)$. The second and third terms are both $O(L^2)$. This completes the proof.

3.2. *Model* 1. We define the hybrid conductance profile

$$\psi_L(t) := \inf_{U \subset \xi_L \,:\, \#(U) \leq t\#(\xi_L)} I_U^3(\xi). \tag{3.14}$$

Note that (3.14) uses the conductance of model 3, but the infimum is over sets $U$ such that $\nu^1(U) \leq t$.

LEMMA 3.3. *The estimate (1.12) of Theorem 1.1 holds for the hybrid profile $\psi_L(t)$.*

PROOF. We repeat the proof of Theorem 1.1. From Lemma 2.2, we know that if $W^3(U) > \gamma|S|$, then $I_U^3 \geq L^{-\varepsilon}$. Let $\gamma, \varepsilon$ be the two constants introduced above. All we then have to prove is the result of Lemma 2.3 adapted to our case. Namely, we need to show that for every $U \subset \xi_L$ such that $\#(U) \leq \frac{1}{2}\#(\xi_L)$ and $W^3(U) \leq \gamma|S|$, we have

$$I_U^3 \geq \delta \min\{L^{-\varepsilon}, \#(U)^{-1/d}\}. \tag{3.15}$$

To prove (3.15), observe that we can again use the bound (2.21). Moreover, thanks to (2.22), we know that $|S| \leq C\#(U)$. Therefore, under the assumption $|S| \leq a|\mathcal{C}_L|$, we have

$$I_U^3 \geq \frac{\kappa\delta(1-a)}{\gamma a} \min\left\{\frac{1}{C\#(U)^{1/d}}, \frac{1}{(\log L)^{d/(d-1)}}\right\}.$$

It remains the case that $|S| > a|\mathcal{C}_L|$. Here in the subcase $|T| \geq 2(1-a)|\mathcal{C}_L|$, things are handled exactly as in Lemma 2.3 and we have the bound $I_U^3 \geq \delta_1(1-a)/\gamma$ as in (2.26). In the subcase $|T| < 2(1-a)|\mathcal{C}_L|$, we use the fact that $W^3(U^c) \geq \#(U^c) \geq \frac{1}{2}\#(\xi_L) \geq \delta|\mathcal{C}_L|$. Therefore, we obtain $W^3(U^c) > \gamma|T|$ for $a$ sufficiently close to 1 [see (2.27)]. Now, observe that by Assumptions (A1), (A2) [see (1.5)] and the bound $W^3(U^c) \geq \frac{1}{2}\#(\xi_L)$ as above, we have the uniform almost sure bound

$$\frac{W^3(U^c)}{W^3(U)} \geq \frac{\#(\xi_L)}{2W^3(\xi_L)} \geq c.$$

In particular, $I_U^3 \geq cI_{U^c}^3$, and the claim follows as in the proof of Lemma 2.3; see (2.28). □

We turn to the proof of Theorem 1.2 for model 1. Using Lemmas 3.1, 3.3 and the arguments in (3.12) and (3.13), it will be sufficient to establish the following lemma.

LEMMA 3.4. *There exists $c > 0$ such that, $\mathcal{P}$-a.s., for $L$ sufficiently large,*

$$(3.16) \qquad \Lambda^1(r) \geq c\psi_L(r)^2, \qquad r \in [\nu_*^1, 1).$$

*Moreover, there exists $C > 0$ such that, $\mathcal{P}$-a.s., for $L$ sufficiently large,*

$$(3.17) \qquad \gamma^1(L) \leq CL^2.$$

PROOF. The proof of (3.16) is similar to the proof of (3.5). We shall use the same notation below. Namely, as in (3.8), we have, for $f \in c_0^+(U)$,

$$\sum_{x,y \in \xi_L} |f(x) - f(y)|Q^1(x,y) = 2\int_0^\infty Q^1(F_t, F_t^c)\,dt.$$

Now, observe that

$$Q^1(F_t, F_t^c) = Q^3(F_t, F_t^c)\frac{W^3(\xi_L)}{\#(\xi_L)} \geq Q^3(F_t, F_t^c) = I_{F_t}^3 \nu^3(F_t).$$

It follows from definition (3.14) that

$$Q^1(F_t, F_t^c) \geq I_{F_t}^3 \nu^3(F_t) \geq \psi_L(\nu^1(F_t))\nu^3(F_t) \geq \psi_L(\nu^1(U))\nu^3(F_t).$$

Therefore, (3.9) becomes

$$\sum_{x,y \in \xi_L} |f(x) - f(y)|Q^1(x,y) \geq 2\psi_L(\nu^1(U))\nu^3(f).$$

Next, repeating the argument after (3.9), we obtain that for any $f \in c_0^+(U)$,

$$(3.18) \quad 2\psi_L(\nu^1(U))\nu^3(f^2) \leq (2\mathcal{E}_1(f))^{1/2}\left(4\sum_{x,y \in \xi_L} f(x)^2 Q^1(x,y)1_{x \neq y}\right)^{1/2}.$$

Using Assumptions (A1), (A2) [see (1.5)] and the definitions (1.6), we have

$$\sum_{y \in \xi_L} Q^1(x,y) 1_{x \neq y} = \frac{w_x^2}{\#(\xi_L)} \leq \frac{w_x^3}{\#(\xi_L)} \leq C \nu^3(x). \tag{3.19}$$

Moreover, there exists $C_1 > 0$ such that for any $x \in \xi_L$, we have

$$\nu^1(x) \leq C_1 \nu^3(x), \qquad x \in \xi_L, \tag{3.20}$$

$\mathcal{P}$-a.s., for $L$ sufficiently large. Indeed, note that $\nu^1(x)/\nu^3(x)$ is given by $\frac{1}{w_x^3} \frac{W^3(\xi_L)}{\#(\xi_L)}$. Since $w_x^3 \geq 1$, the claim (3.20) follows from $\frac{W^3(\xi_L)}{\#(\xi_L)} \leq C_1$, which is a consequence of Assumptions (A1) and (A2). From (3.20), we have $\nu^1(f^2) \leq C_1 \nu^3(f^2)$ for any $f$. In particular, combining (3.18) and (3.19), we have obtained the estimate

$$\psi_L(\nu^1(U)) \leq C \sqrt{\frac{\mathcal{E}_1(f)}{\nu^1(f^2)}} \tag{3.21}$$

for any $f^2 \in c_0^+(U)$. This implies our claim in (3.16).

The proof of (3.17) is a consequence of the comparison estimates

$$\mathrm{Var}_1(f) \leq C \, \mathrm{Var}_3(f), \tag{3.22}$$

$$\mathcal{E}_3(f) \leq C \mathcal{E}_1(f), \tag{3.23}$$

which imply that $\gamma^1(L) \leq C^2 \gamma^3(L)$ [and $\gamma^3(L) \leq CL^2$ follows from (3.6) and Theorem 1.1]. To prove (3.22), we use (3.20):

$$\mathrm{Var}_1(f) \leq \nu^1[(f - \nu^3(f))^2] \leq C_1 \nu^3[(f - \nu^3(f))^2] = C_1 \mathrm{Var}_3(f).$$

The proof of (3.23) is as follows. From Assumptions (A1) and (A2), we have $\frac{\#(\xi_L)}{W^3(\xi_L)} \leq C$. Therefore, for $y \neq x$,

$$\nu^1(x) \mathcal{L}^1(x,y) = \frac{e^{-|x-y|^\alpha}}{\#(\xi_L)} \geq C^{-1} \frac{e^{-|x-y|^\alpha}}{W^3(\xi_L)} = C^{-1} \nu^3(x) \mathcal{L}^3(x,y).$$

This implies that $\mathcal{E}_1(f) \geq C^{-1} \mathcal{E}_3(f)$ for any $f$. □

**4. Proof of Theorem 1.3.** We start with the subdiffusive behavior (case $\alpha > d$ or $\alpha = d$ and $\rho$ small). We show that, $\mathcal{P}_{*,\rho}$-a.s.,

$$\Phi_L^i \leq L^{-1/\delta}, \tag{4.1}$$

for $L$ sufficiently large, where $\Phi_L^i$ is Cheeger's constant and $\delta$ can be taken arbitrarily small (in the case $\alpha = d$, this requires that $\rho$ is accordingly taken to be small). The claims about the Poincaré constant (1.27) and (1.28) will then follow from the simple bound $\gamma^i(L) \geq \frac{1}{2}(\Phi_L^i)^{-1}$, which, in turn, follows from

$$\gamma^i(L) \geq \sup_{U \subset \xi_L : \nu^i(U) \leq \frac{1}{2}} \frac{\operatorname{Var}_i(1_U)}{\mathcal{E}_i(1_U)} = \sup_{U \subset \xi_L : \nu^i(U) \leq \frac{1}{2}} \frac{1 - \nu^i(U)}{I_U^i}.$$

As we see below, the bound (4.1) is a consequence of trapping in isolated regions. For models 1, 3, this can already be achieved at single isolated points. For model 2, we need at least two neighboring points isolated from the rest to produce a subdiffusive behavior. We are going to analyze these situations separately.

4.1. *Case $\alpha > d$. Models* 1, 3. To prove (4.1), we observe that, $\mathcal{P}_{*,\rho}$-a.s., for $L$ sufficiently large, at least one of the $C_W(\log L)^{1/d}$-cubes $W_x \subset \Lambda_L$ has the property that there is exactly one point $x_*$ in $\xi \cap W_x$ such that $d(x_*, W_x^c)$ (Euclidean distance from $x_*$ to the complement of $W_x$) is larger than $\frac{1}{4}C_W(\log L)^{1/d}$. Indeed, let $W$ be a given $C_W(\log L)^{1/d}$-cube and denote by $w$ the unit cube with the same center as $W$ but with volume 1. Let $E$ denote the event that $\xi(w) = 1$ and $\xi(W \setminus w) = 0$. Then,

$$\mathcal{P}_{*,\rho}(E) = \mathcal{P}_{*,\rho}(\xi(w) = 1)\mathcal{P}_{*,\rho}(\xi(W \setminus w) = 0)$$

(4.2) $$\geq \mathcal{P}_{*,\rho}(\xi(w) = 1)\mathcal{P}_{*,\rho}(\xi(W) = 0)$$

(4.3) $$= \rho e^{-\rho} e^{-\rho C_W^d \log L} =: q.$$

Since there are $L^d/C_W^d \log L$-cubes $W_x$ in $\Lambda_L$, the probability that there is no $W_x \subset \Lambda_L$ with the property above is bounded from above by

(4.4) $$(1-q)^{L^d/C_W^d \log L} \leq \exp\left(-\frac{qL^d}{C_W^d \log L}\right).$$

For every $\rho > 0$, we can find sufficiently small $C_W$ such that, for example, $qL^d \geq L^{d/2}$, which implies that (4.4) is summable. Therefore, our claim about the existence of the point $x_*$ follows from the Borel–Cantelli lemma.

Once we have the point $x_*$ as above, we can choose $U = \{x_*\}$. It is simple to check that $\nu^i(U) \leq 1/2$, $\mathcal{P}$-a.s., for $L$ sufficiently large. Hence, $\Phi_L^i \leq I_U^i$. Moreover, since $W^3(U) \geq W^1(U) = 1$, we have $I_U^3 \leq I_U^1$. In conclusion, for $i = 1, 3$,

(4.5) $$\Phi_L^i \leq I_U^i \leq I_U^1 \leq \xi(\Lambda_L) e^{-(1/4C_W(\log L)^{1/d})^\alpha}.$$

Since $\alpha > d$ and $\xi(\Lambda_L) = O(L^d)$ almost surely [see (1.5)], this concludes the proof of (4.1).

4.2. *Case $\alpha > d$. Model 2.* Here, we use the fact that at least one of the $C_W(\log L)^{1/d}$-cubes $W_x \subset \Lambda_L$ has the property that there are exactly two points $x_*, y_*$ in $\xi \cap W_x$ such that both $d(x_*, W_x^c)$ and $d(y_*, W_x^c)$ are larger than $\frac{1}{4}C_W(\log L)^{1/d}$ and such that $d(x_*, y_*)$ is bounded by $\sqrt{d}$. Namely, let $W$ be a given $C_W(\log L)^{1/d}$-cube and denote by $w$ the unit cube with the same center as $W$, as in the argument given above. Let $E'$ denote the event that $\xi(w) = 2$ and $\xi(W \setminus w) = 0$. Then,

$$\begin{aligned}(4.6) \quad \mathcal{P}_{*,\rho}(E') &\geq \mathcal{P}_{*,\rho}(\xi(w) = 2)\mathcal{P}_{*,\rho}(\xi(W) = 0) \\ &= \tfrac{1}{2}\rho^2 e^{-\rho} e^{-\rho C_W^d \log L} =: q'.\end{aligned}$$

For a fixed $\rho$, we can choose $C_W$ sufficiently small in such a way that

$$(4.7) \qquad \exp\left(-\frac{q' L^d}{C_W^d \log L}\right)$$

becomes summable. Therefore, by the same arguments as above, the almost sure existence of the points $x_*, y_*$ follows from the Borel–Cantelli lemma.

Given the points $x_*, y_*$ as above, we observe that, choosing $U = \{x_*, y_*\}$, we have $W^2(U) \geq 2e^{-|x_* - y_*|^\alpha} \geq \delta$ for some $\delta = \delta(\alpha, d)$. However, for $z \in U$ and $z \in U^c$, we have $|z - z'| \geq \frac{1}{4}C_W(\log L)^{1/d}$. It is simple to check that $\nu^2(U) \leq 1/2$, $\mathcal{P}$-a.s., for $L$ large. Hence, we can conclude that

$$(4.8) \qquad \Phi_L^2 \leq I_U^2 \leq \frac{1}{\delta}\xi(\Lambda_L)e^{-(1/4C_W(\log L)^{1/d})^\alpha}.$$

As in (4.5), this implies the subdiffusive estimate (4.1).

4.3. *Case $\alpha = d$ at small density. Models $1, 2, 3$.* We now turn to the case $\alpha = d$. Here, the constant $C_W$ of the $C_W(\log L)^{1/d}$-cubes $W_x$ plays an important role. In this case, we proceed with the same arguments leading to (4.1) in the case $\alpha > d$. Namely, however large the constant $C_W$, using (4.3) and (4.6), we see that if $\rho$ is suitably small [e.g., $\rho < d/(2C_W^d)$], then we can find the desired point $x_*$ or the couple $\{x_*, y_*\}$, as in the cases discussed above, with probability 1. Then, as in (4.5) or (4.8), for any $i = 1, 2, 3$,

$$(4.9) \qquad \Phi_L^i \leq I_U^i \leq C\xi(\Lambda_L)e^{-(C_W^d/4^d)\log L}.$$

Since $\xi(\Lambda_L) = O(L^d)$ almost surely, (4.1) follows by taking $C_W$ sufficiently large (and $\rho$ sufficiently small).

4.4. *Case $\alpha = d$, at high density. Models $1, 2, 3$.* To prove the claim for $\rho$ large, we use the same argument as in the proof of Theorem 1.1. Recall that the only place where the constant $\alpha$ had a role in that proof was in Lemma 2.2. As explained in Remark 3, in the case $\alpha = d$, we need to take $C_W$ sufficiently small so that (2.9) holds with, say, $\varepsilon = \frac{1}{2}$. Also, note that this value of $\varepsilon$ in the estimate for the isoperimetric profile $\varphi_L^i$ in (1.12) is sufficient to prove the desired estimates on mixing times; see (3.13).

Therefore, we need only exploit the fact that if $\rho$ is suitably large, then $C_W$ can be small and we still have that, almost surely, all $C_W(\log L)^{1/d}$-cubes $W_x \subset \Lambda_L$ intersect $\xi_L$; see (2.2). This is possible because if $W_x$ is a $C_W(\log L)^{1/d}$-cube, then

$$\mathcal{P}_{*,\rho}(\xi(W_x) = 0) = e^{-\rho C_W^d \log L}.$$

Therefore, the probability that there exists one such cube with $\xi(W_x) = 0$ and $W_x \cap \Lambda_L \neq \varnothing$ is bounded above by

$$L^d e^{-\rho C_W^d \log L}.$$

Thus, the Borel–Cantelli lemma shows that it suffices to take $\rho > (d + 1)/C_W^d$. This completes the proof of Theorem 1.3.

**5. Examples.** We are going to describe conditions that guarantee the process $\mathcal{P}$ satisfies Assumptions (A1) and (A2). To check the stochastic domination requirement (A1), a very useful criterion is provided by one of the main results of [15], which can be reformulated in our setting as follows. Recall that $B_x, x \in \mathbb{Z}^d$ are the cubes of side $K$ and that $\sigma_x$ is the indicator of the event $\xi(B_x) \geq 1$. For every $x \in \mathbb{Z}^d$ and $D > 0$, let $U_{D,x}$ denote the set $U_{D,x} = \{y \in \mathbb{Z}^d : |y - x| > D\}$.

LEMMA 5.1. *Suppose that there exist $D > 0$ such that for all $x \in \mathbb{Z}^d$,*

(5.1) $$\mathcal{P}(\sigma_x = 1 | \sigma = \zeta \text{ on } U_{D,x}) \geq p$$

*for $\mathcal{P}$-a.a. $\zeta \in \{0,1\}^{\mathbb{Z}^d}$, where $p = p(D,K)$ is such that $\lim_{K \to \infty} p(D,K) = 1$. Then there exists $\rho = \rho(D,K)$ with $\lim_{K \to \infty} \rho(D,K) = 1$, such that the random field $\{\sigma_x\}$ stochastically dominates the product Bernoulli process with parameter $\rho$.*

PROOF. This is a special case of Theorem 1.3 in [15]. □

Concerning Assumption (A2), we know that it is satisfied by any process $\mathcal{P}$ such that $\mathcal{P} \preceq \mathcal{P}_{*,\rho}$ for some $\rho > 0$. This is proved in Appendix A, Proposition A.1. More generally, we expect that Assumption (A2) is satisfied by any $\mathcal{P}$ which is stochastically dominated by a process $\widetilde{\mathcal{P}}$ with good mixing properties, for example, exponential tree decay of correlations. We turn to some specific examples.

5.1. *Poisson processes.* Suppose that $\mathcal{P} = \mathcal{P}_{*,\rho}$ is the homogeneous Poisson point process with intensity $\rho > 0$. Assumption (A1) then is obviously satisfied since $\xi(B_x)$ are i.i.d. Poisson random variables with parameter $\rho K^d$ and $\sigma_x$ are independent Bernoulli variables with $p = 1 - e^{-\rho K^d}$. It follows that Assumption (A1) is satisfied by any $\mathcal{P}$ such that $\mathcal{P} \succeq \mathcal{P}_{*,\rho}$ for some $\rho > 0$. Moreover, Assumption (A2) is satisfied by any $\mathcal{P}$ such that $\mathcal{P} \preceq \mathcal{P}_{*,\rho}$ for some $\rho > 0$. In particular, if $\mathcal{P}$ is any process such that

$$\mathcal{P}_{*,\rho_1} \preceq \mathcal{P} \preceq \mathcal{P}_{*,\rho_2} \tag{5.2}$$

with some $0 < \rho_1 < \rho_2 < \infty$ then both Assumptions (A1) and (A2) hold. The domination (5.2), holds in particular, for nonhomogeneous Poisson processes with intensity function $\varphi(x)$ such that $\rho_1 \leq \varphi(x) \leq \rho_2$ (see, e.g., [12]).

5.2. *Thinning of point processes with uniform bounds on the local density.* Consider a point process $\xi$ such that, $\mathcal{P}$-a.s.,

$$1 \leq \xi(\ell x + \Lambda_\ell) \leq n \qquad \forall x \in \mathbb{Z}^d, \tag{5.3}$$

for suitable constants $n, \ell > 0$, where $\Lambda_\ell = [-\frac{\ell}{2}, \frac{\ell}{2}]^d$. Given $p \in (0,1]$, let $\hat{\xi}$ be the $p$-thinning of $\xi$, that is, $\hat{\xi}$ is obtained from $\xi$ by erasing points of $\xi$ independently with probability $1 - p$. Note that $\xi = \hat{\xi}$ if $p = 1$. The process $\hat{\xi}$ can model both crystal/quasicrystal structures ($p = 1$) and their variants due to defects [$p \in (0,1)$]. Trivially, $\hat{\xi}$ satisfies Assumptions (A1), (A2). A typical example of point process $\hat{\xi}$ is given by the diluted $\mathbb{Z}^d$, defined as the $p$-thinning of $\xi \equiv \mathbb{Z}^d$.

5.3. *High-temperature/low-fugacity gas.* Consider a Gibbsian random point field described by the formal Hamiltonian function

$$H(\xi) = \sum_{x,y \in \xi} \varphi(x - y), \tag{5.4}$$

where $\varphi : \mathbb{R}^d \to \mathbb{R}$ is an even function (the two-body potential). It is known that under suitable hypothesis on $\varphi$ and for sufficiently small values of the inverse temperature $\beta$ and of the fugacity $\lambda$, one can apply cluster expansion techniques to obtain a well-defined Gibbs field $\mathcal{P}_{\beta,\lambda}$ in the usual DLR sense [24]. We now consider the case of nonnegative finite range potentials in detail. We comment briefly on other models afterward.

5.3.1. *Nonnegative, finite range potential.* Suppose that $\varphi : \mathbb{R}^d \to \mathbb{R}$ is a measurable even function such that $\varphi \geq 0$ and $\varphi(x) = 0$ for $|x| > R$, for some $R < \infty$. A uniformly convergent cluster expansion for such functions that has been obtained by several authors. In particular, at sufficiently small values of $\beta, \lambda$, there exists a unique Gibbs measure $\mathcal{P} = \mathcal{P}_{\beta,\lambda}$ for the interaction (5.4).

In [26], this is derived, together with exponential clustering properties for the random field $\mathcal{P}$ that hold uniformly in the boundary conditions outside a given region. We write $\mathcal{P}_\Lambda^\eta$ for the Gibbs measure in a bounded Borel subset $\Lambda \subset \mathbb{R}^d$ with boundary condition $\eta$ as follows. Let $\Omega_\Lambda$ denote the set of finite subsets of $\Lambda$, endowed with the $\sigma$-algebra $\mathcal{F}_\Lambda$ generated by the counting functions $N_A : \xi \to \#(\xi \cap A)$, $A \subset \Lambda$. Then, if $f$ is a measurable function on $\Omega_\Lambda$, we define

$$(5.5) \qquad \mathbb{E}_{\mathcal{P}_\Lambda^\eta}[f] = \frac{1}{Z_\Lambda^\eta} \sum_{n=0}^\infty \frac{\lambda^n}{n!} \int_{\Lambda^n} e^{-\beta H_\Lambda^\eta(\omega)} f(\omega) \, d\omega,$$

where $f$ has been identified with a symmetric function on $\bigcup_{n \geq 0} \Lambda^n$, $Z_\Lambda^\eta$ is the normalizing constant and, for any finite $\omega \subset \Lambda$ and any locally finite $\eta \subset \mathbb{R}^d$,

$$(5.6) \qquad H_\Lambda^\eta(\omega) = \sum_{\substack{\{x,y\} \subset \omega \cup (\eta \cap \Lambda^c) : \\ \{x,y\} \cap \Lambda \neq \varnothing}} \varphi(x-y).$$

The following estimates have been established in [3], Corollary 2.4 and Corollary 2.5, based on the expansion presented in [26]. We give some preliminary notation. The *support* of $f$ is the smallest $\Lambda$ such that $f$ is $\mathcal{F}_\Lambda$ measurable and is denoted by $\Lambda_f$. Moreover, $\bar\Lambda_f$ stands for its Euclidean enlargement by $R$, where $R$ is the range of the interaction, that is, $\bar\Lambda_f = \{x \in \mathbb{R}^d : d(x, \Lambda_f) \leq R\}$. For any $\eta, \tau \in \Omega$, we set $\eta \Delta \tau = (\eta \cup \tau) \setminus (\eta \cap \tau)$ for the symmetric difference. Corollary 2.5 of [3] is then stated as follows.

LEMMA 5.2. *Let $\beta, \lambda$ be such that*

$$(5.7) \qquad \frac{\lambda \epsilon(\beta)}{1 - 2\lambda \epsilon(\beta)} < 1, \qquad \epsilon(\beta) := e \int_{\mathbb{R}^d} (1 - e^{-\beta \varphi(x)}) \, dx.$$

*There then exist $C < \infty$ and $m > 0$ such that for any bounded $\Lambda \subset \mathbb{R}^d$, any local function $f$ and any pair of boundary conditions $\eta, \tau$ satisfying $d(\Lambda_f, \eta \Delta \tau) > 3R$ and $|\bar\Lambda_f| \leq \exp[m(d(\Lambda_f, \eta \Delta \tau) - R)]$, we have*

$$(5.8) \qquad |\mathbb{E}_{\mathcal{P}_\Lambda^\eta}[f] - \mathbb{E}_{\mathcal{P}_\Lambda^\tau}[f]| \leq C \left( \sup_{\eta'} \mathbb{E}_{\mathcal{P}_\Lambda^{\eta'}}[|f|] \right) e^{-m d(\Lambda_f, \eta \Delta \tau)}.$$

The above result implies, in particular, that for $\beta, \lambda$ satisfying (5.7), the Gibbs field $\mathcal{P} = \mathcal{P}_{\beta, \lambda}$ is unique. Moreover, $\mathcal{P}_{\beta, \lambda}$ is stationary and ergodic.

LEMMA 5.3. *Assume (5.7). The Gibbs measure $\mathcal{P} = \mathcal{P}_{\beta, \lambda}$ then satisfies Assumptions* (A1) *and* (A2).

PROOF. We start with Assumption (A2). This follows from the stochastic domination $\mathcal{P}_{\beta,\lambda} \preceq \mathcal{P}_{*,\lambda}$, which is a consequence of nonnegativity of $\varphi$ (repulsive interaction); see [12].

We turn to the proof of Assumption (A1). We shall establish (5.1) for the Gibbs measure $\mathcal{P}$. Let us first observe that from ergodicity, it follows that $\mathcal{P}(\sigma_x = 1) \to 1$ as $K \to \infty$. We then use Lemma 5.2 to obtain the desired estimate in (5.1). Note that if $|x - y| > D$ with $D$ sufficiently large, then $d(B_x, B_y) > DK/2$ for all $K$. Let $\Lambda = \Lambda(D, K)$ denote the set

$$\Lambda = \{z \in \mathbb{R}^d : d(z, B_x) \leq DK/2\}.$$

We write

(5.9) $\quad \mathcal{P}(\sigma_x = 1 | \sigma = \zeta \text{ on } U_{D,x}) = \int \mathbb{E}_{\mathcal{P}_\Lambda^\eta}[\sigma_x] \mathcal{P}(d\eta | \sigma = \zeta \text{ on } U_{D,x}).$

From (5.8), taking $f = \sigma_x$, $\Lambda_f = B_x$ and $\Lambda = \Lambda(D, K)$, we see that for $K$ sufficiently large, $\mathbb{E}_{\mathcal{P}_\Lambda^\eta}[\sigma_x] \geq \mathbb{E}_{\mathcal{P}_\Lambda^\tau}[\sigma_x] - Ce^{-mDK/2}$ for any pair of boundary conditions $\eta, \tau$, with some independent constant $C$. It then follows that

$$\mathbb{E}_{\mathcal{P}_\Lambda^\eta}[\sigma_x] \geq \mathbb{E}_\mathcal{P}[\sigma_x] - Ce^{-mDK/2}$$

uniformly in $\eta$. Therefore, using (5.9), we obtain

(5.10) $\quad \mathcal{P}(\sigma_x = 1 | \sigma = \zeta \text{ on } U_{D,x}) \geq \mathbb{E}_\mathcal{P}[\sigma_x] - Ce^{-mDK/2}.$

The claim now follows from $\mathbb{E}_\mathcal{P}[\sigma_x] = \mathcal{P}(\sigma_x = 1) \to 1$ as $K \to \infty$. This completes the proof of Assumption (A1). $\square$

REMARK 6. All of the above examples fulfill Assumptions (A1) and (A2) for every $\alpha > 0$. Hence, due to Remark 1, there exists $C > 0$ such that $\gamma^1(L) \geq CL^2$, $\mathcal{P}$-a.s., for $L$ sufficiently large.

5.4. *Other examples.* We expect Assumptions (A1) and (A2) to hold for Gibbsian fields $\mathcal{P}_{\beta,\lambda}$ whenever one has a uniformly convergent high-temperature/low-fugacity expansion with clustering properties that hold uniformly in the boundary conditions, as in Lemma 5.2. The latter is known to be the case for some models with multibody interactions under the assumption that the pair potential $\varphi$ satisfies $\varphi(x) = +\infty$ for $|x| < R_0$, for some $R_0 < \infty$ (hard-core interactions) and under some mild additional assumptions [23] (in particular, one can remove the positivity and finite range requirement on $\varphi$). For more general models with only pair interaction, such as the one considered in [14], where $\varphi$ is only assumed to be stable and exponentially decaying at infinity, the clustering property derived in Theorem 2 of [14] is not sufficient to establish Assumptions (A1), (A2) here because of the lack of uniformity in the boundary condition. In particular, a uniform result, as in Lemma 5.2, is not available in this case.

Finally, it would be interesting to investigate the validity of Assumptions (A1) and (A2) in other classes of point processes. The class of *determinantal processes* received much attention recently (see, e.g., [1] and references therein). Due to the negative association property of these processes (repulsion), the most delicate issue here seems to be Assumption (A1). We are not aware of results in this direction for determinantal point processes in the continuum. On the other hand, for lattice determinantal processes, simple explicit criteria are known for stochastic domination from above and from below in terms of i.i.d. Bernoulli random variables [17].

5.5. *Palm distribution.* In [10], the authors consider a random walk on the support of a marked simple point process whose jump rates decay exponentially in the jump length and depend via a Boltzmann-type factor on the (energy) marks. The law of the process is the Palm distribution associated with a stationary ergodic marked simple point process. Since the Boltzmann-type factor in the marks is bounded from above and below by positive constants, the estimates of Cheeger's constant, spectral gap and mixing time for the random walk confined in a cubic box of side $L$ reduce to the case of zero energy marks and are hence covered by the following discussion.

Recall that $\Omega$ denotes the Borel space of locally finite subsets $\xi \subset \mathbb{R}^d$, endowed with the $\sigma$-algebra $\mathcal{F}$ generated by the counting variables $N_\Lambda(\xi) = \#(\xi \cap \Lambda)$, and define $\Omega_0$ as the Borel subset of $\Omega$ given by the subsets $\xi$ containing the origin. Given a stationary simple point process on $\mathbb{R}^d$ with law $\mathcal{P}$ and finite intensity $\rho$, that is,

$$\rho := \mathbb{E}_\mathcal{P}(\xi(\Lambda_1)) < \infty,$$

the associated Palm distribution $\mathcal{P}_0$ is the probability measure on $\Omega_0$ such that

$$(5.11) \qquad \mathcal{P}_0(A) = \frac{1}{\rho} \int_\Omega \mathcal{P}(d\xi) \sum_{z \in \xi_1} \chi_A(\tau_z \xi) \qquad \forall A \subset \Omega_0 \text{ Borel},$$

where $\Lambda_1 = [-\frac{1}{2}, \frac{1}{2}]^d$, $\xi_1 = \xi \cap \Lambda_1$, $\chi_A$ denotes the characteristic function of the event $A$ and $\tau_z \xi$ denotes the translated subset $\xi - z$.

One can directly apply the results described in the Introduction to the case $\mathcal{P}_0$. Since it is usually simpler to deal with the original law $\mathcal{P}$ than with the associated Palm distribution $\mathcal{P}_0$, it is useful to obtain results for $\mathcal{P}_0$ under suitable assumptions on $\mathcal{P}$ instead of $\mathcal{P}_0$. We prove that if $\mathcal{P}$ satisfies Assumptions (A1) and (A2), then the conclusions of Theorems 1.1 and 1.2 still hold for $\mathcal{P}_0$.

LEMMA 5.4. *Suppose that $\mathcal{P}$ satisfies Assumptions (A1) and (A2) and that $\alpha < d$. For every $\varepsilon > 0$, there exists $\delta > 0$ such that for all $i = 1, 2, 3$, $\mathcal{P}_0$-a.s.,*

$$\varphi_L^i(t) \geq \delta \min\left\{\frac{1}{L^\varepsilon}, \frac{1}{t^{1/d} L}\right\}, \qquad 0 < t \leq \frac{1}{2}, \tag{5.12}$$

*for all $L$ sufficiently large. In particular, there exists $C < \infty$ such that, $\mathcal{P}_0$-a.s.,*

$$\Phi_L^i \geq \delta L^{-1}, \tag{5.13}$$

$$\tau^i(L) \leq C L^2 \tag{5.14}$$

*for all $L$ sufficiently large.*

PROOF. As discussed in Section 3, (5.13) and (5.14) are a consequence of (5.12). In order to prove (5.12), consider the event $A_{\delta,\varepsilon,L_0}$ given by the subsets $\xi \in \Omega$ satisfying (5.12) for $L \geq L_0$, $L \in \mathbb{N}$. Due to Remark 4, given $\varepsilon > 0$, there exists $\delta > 0$ such that

$$\mathcal{P}(\exists L_0 > 0 : \tau_z \xi \in A_{\delta,\varepsilon,L_0} \text{ for all } z \in \xi_1) = 1.$$

In particular, due to (5.11), $\lim_{L_0 \uparrow \infty} \mathcal{P}_0(A_{\delta,\varepsilon,L_0}) = 1$, thus implying (5.12) for $L$ sufficiently large $\mathcal{P}_0$-a.s. $\square$

## APPENDIX A: ON THE VALIDITY OF ASSUMPTION (A2)

Recall that, for any bounded Borel set $A \subset \mathbb{R}^d$,

$$R_A(\xi) = \sum_{x \in \xi \cap A} \sum_{y \in \xi} e^{-|x-y|^\alpha}. \tag{A.1}$$

Also, recall that $\Lambda_\ell = [-\frac{\ell}{2}, \frac{\ell}{2}]^d$, $\ell \in \mathbb{N}$, and that $\mathcal{P}_{*,\rho}$ denotes the homogeneous Poisson point process with intensity $\rho > 0$.

PROPOSITION A.1. *There exists a constant $c_1 = c_1(\rho, \alpha, d)$ such that for every $n \in \mathbb{N}$,*

$$\mathcal{P}_{*,\rho}(R_{\Lambda_\ell} \geq \gamma \ell^d) \leq C_n \ell^{-nd}, \qquad \forall \gamma > c_1, \ell \in \mathbb{N}, \tag{A.2}$$

*where $C_n$ is a finite constant depending on $n$ and $\rho, \alpha, d$. In particular, Assumption (A2) is satisfied by any $\mathcal{P}$ such that $\mathcal{P} \preceq \mathcal{P}_{*,\rho}$ for some $\rho > 0$.*

PROOF. Let us first show that Assumption (A2) is satisfied by any $\mathcal{P} \preceq \mathcal{P}_{*,\rho}$ once we have the bound (A.2). For $\varepsilon > 0$, recall the definition of the cubes $V_x = L^\varepsilon x + [0, L^\varepsilon)^d$, $x \in \mathbb{Z}^d$. Since $R_{V_x}$ is an increasing function, it is sufficient to establish the bound (1.4) for $\mathcal{P}_{*,\rho}$. This, in turn, is an immediate consequence of (A.2) since a union bound shows that the probability that at least one of the $V_x$ in the partition of $\Lambda_L$ is such that $R_{V_x} \geq CL^{\varepsilon d}$ for a sufficiently large constant $C$ is bounded from above by

$$L^d C_n L^{-\varepsilon dn}.$$

Therefore, we can choose, for example, $n = n(\varepsilon) = \frac{d+2}{\varepsilon d}$ to obtain that the probability of violating Assumption (A2) is summable in $L \in \mathbb{N}$ and the claim follows from the Borel–Cantelli lemma.

We turn to the proof of (A.2). To simplify the notation, we shall use the convention that $c, c', c'', \ldots$ stand for positive constants depending only on the parameters $\rho, \alpha, d$ and $n$ (but not on $\ell$), whose values can change from line to line.

We define

$$S_\ell(\xi) = \sum_{u \in \Lambda_\ell \cap \mathbb{Z}^d} \sum_{v \in \mathbb{Z}^d} e^{-|u-v|^\alpha} \xi(Q_u) \xi(Q_v),$$

where $Q_u := u + [-1/2, 1/2]^d$, $u \in \mathbb{Z}^d$. It is simple to check that, for $c = c(\alpha, d)$,

$$R_{\Lambda_\ell}(\xi) \leq c S_\ell(\xi).$$

Hence, it is enough to prove Proposition A.1 with $R_{\Lambda_\ell}$ replaced by $S_\ell$.

Note that $\mathbb{E}_{*,\rho}[S_\ell] \leq c_1 \ell^d$ for some $c_1$ depending on $\rho, \alpha, d$ (here and below, $\mathbb{E}_{*,\rho}$ stands for expectation w.r.t. $\mathcal{P}_{*,\rho}$). Therefore, for $\gamma$ sufficiently large, we can write

$$\begin{aligned}(\text{A.3})\qquad \mathcal{P}_{*,\rho}(S_\ell > \gamma \ell^d) &\leq \mathcal{P}_{*,\rho}(S_\ell - \mathbb{E}_{*,\rho}[S_\ell] > (\gamma/2)\ell^d) \\ &\leq (\gamma/2)^{-2n} \ell^{-2nd} \mathbb{E}_{*,\rho}[(S_\ell - \mathbb{E}_{*,\rho}[S_\ell])^{2n}].\end{aligned}$$

Thus, it will be sufficient to show that

(A.4) $$\mathbb{E}_{*,\rho}[(S_\ell - \mathbb{E}_{*,\rho}[S_\ell])^{2n}] \leq c \ell^{nd}.$$

Let us define

$$F_{u,v}(\xi) = \xi(Q_u)\xi(Q_v) - \mathbb{E}_{*,\rho}[\xi(Q_u)\xi(Q_v)].$$

We can then write

$$\begin{aligned}\mathbb{E}_{*,\rho}&[(S_\ell - \mathbb{E}_{*,\rho}[S_\ell])^{2n}] \\ &= \sum{}^* e^{-|u_1-v_1|^\alpha - |u_2-v_2|^\alpha - \cdots - |u_{2n}-v_{2n}|^\alpha} \\ &\qquad \times \mathbb{E}_{*,\rho}[F_{u_1,v_1}(\xi) F_{u_2,v_2}(\xi) \cdots F_{u_{2n},v_{2n}}(\xi)],\end{aligned}$$

where the sum $\sum^*$ is defined as

$$\sum{}^* = \sum_{u_1 \in \Lambda_\ell \cap \mathbb{Z}^d} \sum_{u_2 \in \Lambda_\ell \cap \mathbb{Z}^d} \cdots \sum_{u_{2n} \in \Lambda_\ell \cap \mathbb{Z}^d} \sum_{v_1 \in \mathbb{Z}^d} \sum_{v_2 \in \mathbb{Z}^d} \cdots \sum_{v_{2n} \in \mathbb{Z}^d}.$$

Note that if
$$\mathbb{E}_{*,\rho}[F_{u_1,v_1}(\xi)F_{u_2,v_2}(\xi)\cdots F_{u_{2n},v_{2n}}(\xi)]$$
is not zero, then for each $i \in \{1,\ldots,2n\}$, there exists $j \in \{1,\ldots,2n\}$ with $i \neq j$ such that $\{u_i,v_i\} \cap \{u_j,v_j\} \neq \varnothing$. Indeed, if there is an isolated pair $u_i, v_i$, then the expression vanishes by independence. Hence, using the simple bound [for some $c = c(n,\rho)$]
$$|\mathbb{E}_{*,\rho}[F_{u_1,v_1}(\xi)F_{u_2,v_2}(\xi)\cdots F_{u_{2n},v_{2n}}(\xi)]| \leq c$$
$$\forall u_1,u_2,\ldots,u_{2n},v_1,v_2,\ldots,v_{2n},$$
if we write $\sum^{**}$ for the sum $\sum^{*}$ restricted to $u_1,\ldots,u_{2n},v_1,\ldots,v_{2n}$ such that the property mentioned above is satisfied, we obtain the bound
$$\mathbb{E}_{*,\rho}[(S_\ell - \mathbb{E}_{*,\rho}[S_\ell])^{2n}] \leq c\sum\nolimits^{**} e^{-|u_1-v_1|^\alpha - |u_2-v_2|^\alpha - \cdots - |u_{2n}-v_{2n}|^\alpha}.$$
In order to complete the proof, we need to show that

(A.5) $$\sum\nolimits^{**} e^{-|u_1-v_1|^\alpha - |u_2-v_2|^\alpha - \cdots - |u_{2n}-v_{2n}|^\alpha} \leq c\ell^{nd}.$$

Let us first observe that the contribution to the left-hand side of (A.5) of addenda such that $\|v_i\|_\infty > \ell$ for some $i$ is bounded by some constant $c$. In fact, suppose, for example, that $\|v_1\|_\infty > \ell$. Since $\|u_1\|_\infty \leq \ell/2$, we can bound
$$|u_1 - v_1| \geq c\|u_1 - v_1\|_\infty \geq c(\|v_1\|_\infty - \|u_1\|_\infty) \geq c'\ell + c''\|v_1\|_\infty.$$
Hence,
$$\sum\nolimits^{**} \chi(\|v_1\|_\infty > \ell) e^{-|u_1-v_1|^\alpha - |u_2-v_2|^\alpha - \cdots - |u_{2n}-v_{2n}|^\alpha}$$
$$\leq e^{-c'\ell^\alpha} \sum_{u_1 \in \Lambda_\ell \cap \mathbb{Z}^d} \sum_{v_1 \in \mathbb{Z}^d} e^{-c''\|v_1\|_\infty^\alpha} \left(\sum_{u \in \Lambda_\ell \cap \mathbb{Z}^d} \sum_{v \in \mathbb{Z}^d} e^{-|u-v|^\alpha}\right)^{2n-1}$$
$$\leq c e^{-c'\ell^\alpha} \ell^{2dn} \leq c'''.$$

Using this observation, it is enough to prove (A.5) when the sum $\sum^{**}$ is restricted to $v_1,\ldots,v_{2n}$ in $[-\ell,\ell]^d \cap \mathbb{Z}^d$. Of course, we may also extend to $u_1,\ldots,u_{2n} \in [-\ell,\ell]^d \cap \mathbb{Z}^d$. Hence, Proposition A.1 follows from Lemma A.2 below. □

LEMMA A.2. *Let $\sum^\star$ denote the sum over $u_1,u_2,\ldots,u_{2n},v_1,v_2,\ldots,v_{2n}$ in $[-\ell,\ell]^d \cap \mathbb{Z}^d$ such that for each $i \in \{1,\ldots,2n\}$, there exists $j \in \{1,\ldots,2n\}$ with $i \neq j$ and $\{u_i,v_i\} \cap \{u_j,v_j\} \neq \varnothing$. Then,*

(A.6) $$\sum\nolimits^\star e^{-|u_1-v_1|^\alpha - |u_2-v_2|^\alpha - \cdots - |u_{2n}-v_{2n}|^\alpha} \leq c\ell^{nd}.$$

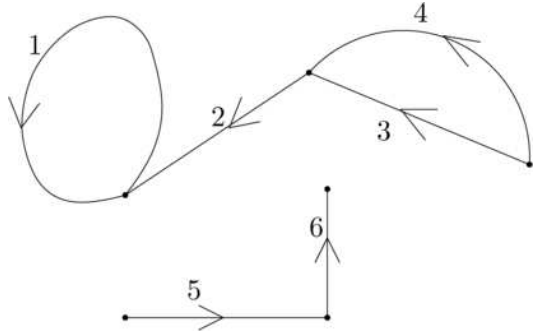

FIG. 4. *Example of a graph $\mathcal{G} \in \Theta_{\{1,2,\ldots,6\}}$.*

PROOF. The proof is based on a combinatorial argument and it is convenient to start by explaining the needed graph-theoretic notation. Given $W \subset \mathbb{N}_+$, we denote by $\Theta_W$ the family of oriented graphs with $|W|$ edges ($|W|$ is the cardinality of the set $W$) such that:

(1) each edge is oriented and labeled by a number $i \in W$ (different edges are labeled by different numbers);

(2) each connected component contains at least two edges.

See Figure 4 for an example. We shall take $W \subset \{1, \ldots, 2n\}$ and will use a graph in $\Theta_W$ to describe the dependence structure between points $(u_1, v_1), \ldots, (u_{2n}, v_{2n})$, as explained below. Graphs in $\Theta_W$ are thought of up to isomorphisms, that is, up to bijective maps from the vertex set of one graph to the vertex set of the other that conserve the orientation and labeling of the edges.

To a graph $\mathcal{G} \in \Theta_W$, with $W \subset \{1, 2, \ldots, 2n\}$, we associate the family $\Omega_\mathcal{G}$ of labeled oriented graphs $G$ satisfying the following properties:

(1) $G$ is isomorphic to $\mathcal{G}$;

(2) $G$ has vertex set $\{u_i\}_{i \in W} \cup \{v_i\}_{i \in W}$, where $u_i, v_i \in [-\ell, \ell]^d \cap \mathbb{Z}^d$;

(3) for each $i \in W$, $u_i$ is connected to $v_i$ by the oriented edge (from $u_i$ to $v_i$) labeled by $i$.

With these definitions, it becomes clear that terms in the sum $\sum^\star$ defined in Lemma A.2 can be enumerated by first enumerating all graphs $\mathcal{G} \in \Theta_W$, with $W = \{1, \ldots, 2n\}$, and then enumerating all graphs $G \in \Omega_\mathcal{G}$ corresponding to that $\mathcal{G}$. For instance, if $n = 3$, the choice of the graph $\mathcal{G}$ given in Figure 4 corresponds to the constraints $u_1 = v_1 = v_2 = a_1, u_2 = v_3 = v_4 = a_2, u_3 = u_4 = a_3, u_5 = a_4, v_5 = u_6 = a_5, v_6 = a_6$, with arbitrary distinct points $a_1, \ldots, a_6 \in [-\ell, \ell]^d \cap \mathbb{Z}^d$. In general, the fact that $\mathcal{G}$ and $G$ are isomorphic implies that $u_{i_1}, u_{i_2}, \ldots, u_{i_k}$ and $v_{j_1}, v_{j_2}, \ldots, v_{j_s}$ must all coincide whenever there exists a vertex $a$ of the graph $\mathcal{G}$ such that the oriented edges labeled by $i_1, \ldots, i_k$ exit from $a$ and the oriented edges labeled by $j_1, \ldots, j_s$ enter in *a.* Given $G \in \Omega_{\mathcal{G}}$, we define

(A.7) $$F(G) = \prod_{i \in W} e^{-|u_i - v_i|^\alpha}.$$

In particular, in the case of Figure 4, we would obtain

$$\sum_{G \in \Omega_{\mathcal{G}}} F(G) = \sum_{\substack{a_1, a_2, \ldots, a_6 \in [-\ell, \ell]^d \cap \mathbb{Z}^d: \\ a_i \text{ all distinct}}} \prod_{i=1}^{6} e^{-|u_i - v_i|^\alpha} \chi(u_1 = v_1 = v_2 = a_1,$$
$$u_2 = v_3 = v_4 = a_2,$$
$$u_3 = u_4 = a_3, u_5 = a_4,$$
$$v_5 = u_6 = a_5, v_6 = a_6).$$

In general, thanks to the above definitions, we can write the sum in Lemma A.2 as

(A.8) $$\sideset{}{^\star}\sum e^{-|u_1-v_1|^\alpha - |u_2-v_2|^\alpha - \cdots - |u_{2n}-v_{2n}|^\alpha} = \sum_{\mathcal{G} \in \Theta_{\{1,2,\ldots,2n\}}} \sum_{G \in \Omega_{\mathcal{G}}} F(G).$$

Next, note that if $\mathcal{G}_1, \mathcal{G}_2, \ldots, \mathcal{G}_k$ are the connected components of a given $\mathcal{G}$, then

(A.9) $$\sum_{G \in \Omega_{\mathcal{G}}} F(G) = \prod_{i=1}^{k} \left( \sum_{G \in \Omega_{\mathcal{G}_i}} F(G) \right).$$

By definition, each connected component of any graph $\mathcal{G} \in \Theta_{\{1,2,\ldots,2n\}}$ must contain at least two edges. Since $\mathcal{G}$ in (A.8) has $2n$ edges, the graph $\mathcal{G}$ has at most $n$ connected components. Thus, the number $k$ in (A.9) is not larger than $n$. Also, note that the number of graphs $\mathcal{G} \in \Theta_{\{1,2,\ldots,2n\}}$ only depends on $n$. Thanks to (A.8) and these last observations, in order to complete the proof, it is enough to prove that

(A.10) $$\sum_{G \in \Omega_{\mathcal{G}}} F(G) \leq c\ell^d$$

for each connected graph $\mathcal{G} \in \Theta_W$, for any $W \subset \{1, 2, \ldots, 2n\}$.

Fix a connected graph $\mathcal{G} \in \Theta_W$ and let $m$ denote the number of its vertices. Note that if $m = 1$, then (A.10) follows immediately. Thus, we will assume that $m \geq 2$. Denote by $\mathcal{G}'$ an arbitrary spanning tree of $\mathcal{G}$, that is, a connected subgraph of $\mathcal{G}$ with the same vertex set which has no cycles. Let $\Omega_{\mathcal{G}'}$ denote the corresponding family of graphs, as defined above. Since $\mathcal{G}'$ is obtained from $\mathcal{G}$ by removing some edges (if necessary), from the definition (A.7) we have

$$\sum_{G \in \Omega_{\mathcal{G}}} F(G) \leq \sum_{G \in \Omega_{\mathcal{G}'}} F(G).$$

Hence, it is enough to prove (A.10) when $\mathcal{G}$ is a tree.

We prove this statement by induction over the number $m$ of vertices of the tree $\mathcal{G}$. If $m = 2$, the statement is straightforward since

$$\sum_{G \in \Omega_\mathcal{G}} F(G) = \sum_{u_1, v_1 \in [-\ell, \ell]^d \cap \mathbb{Z}^d \, : \, u_1 \neq v_1} e^{-|u_1 - v_1|^\alpha} = O(\ell^d).$$

Suppose, now, that $m > 2$. Take a leaf $a$ of $\mathcal{G}$, that is, a vertex that is connected to only one other vertex $b$ by an edge $e = (b, a)$ or $e = (a, b)$. Consider the new tree $\hat{\mathcal{G}}$ obtained from $\mathcal{G}$ by removing the vertex $a$ and the edge $e$. Then,

$$\sum_{G \in \Omega_\mathcal{G}} F(G) \leq \sum_{\hat{G} \in \Omega_{\hat{\mathcal{G}}}} F(\hat{G}) \left( \sum_{z \in [-\ell, \ell]^d \cap \mathbb{Z}^d} e^{-|w(\hat{G}) - z|^\alpha} \right),$$

where $w(\hat{G})$ is the vertex of $\hat{G}$ corresponding to the vertex $b$ in $\hat{\mathcal{G}}$ via the isomorphism between $\hat{\mathcal{G}}$ and $\hat{G}$. The last factor in the expression above is bounded by $\sum_{z \in \mathbb{Z}^d} e^{-|z|^\alpha} < \infty$. Therefore,

$$\sum_{G \in \Omega_\mathcal{G}} F(G) \leq c \sum_{\hat{G} \in \Omega_{\hat{\mathcal{G}}}} F(\hat{G}).$$

Since the number of edges in $\mathcal{G}$ is at most $2n$, we can iterate this estimate down to the case $m = 2$ and the proof is complete. □

## APPENDIX B: PERCOLATION RESULTS

In this appendix, we prove some properties concerning the maximal open cluster in a finite box for site percolation with parameter $p$ close to 1. Similar results hold for each supercritical $p$; see [22]. For the reader's convenience, we give a simple and essentially self-contained proof for the case of large $p$.

We consider site percolation on $\mathbb{Z}^d$, $d \geq 2$, that is, we have i.i.d. Bernoulli random variables $\omega(x), x \in \mathbb{Z}^d$, with parameter $p \in (0, 1)$. As usual, a point $x$ is said to be *open* if $\omega(x) = 1$, *closed* if $\omega(x) = 0$. We denote by $B_n$ the box $\{1, 2, \ldots, n\}^d$ and consider the natural graph structure inherited by $\mathbb{Z}^d$, that is, two points $x, y \in B_n$ are joined by an edge iff $|x - y| = 1$, where $|x - y|$ denotes Euclidean distance. A set $A \subset B_n$ will be called $B_n$-*connected* if it is connected with respect to this structure. Moreover, a set $A \subset B_n$ is called 2-*connected* if for every $x, y \in A$, there exists a path $x = z_1, \ldots, z_m = y$ such that $z_i \in A$ and $|z_i - z_{i+1}| \leq 2$ for every $i = 1, \ldots, m$. For $A \subset B_n$, we define $d_\infty(A) = \max\{\|x - y\|_\infty, x, y \in A\}$, the diameter in the $\ell_\infty$-norm. The $B_n$-connected components of the set of open vertices in $B_n$ are called *open clusters* or simply *clusters*. We introduce the events $\mathcal{A}_n$, $\mathcal{B}_n$ and $\mathcal{C}_n(\kappa)$ as follows: $\mathcal{A}_n$ is the event that there exists at most one open cluster $C$ in $B_n$ such that $d_\infty(C) \geq [n/10]$; $\mathcal{B}_n$ is the event that there exists an open crossing cluster in $B_n$, that is, an open cluster intersecting all of the faces of $B_n$; given $\kappa \in (0,1)$, $\mathcal{C}_n(\kappa)$ denotes the event that there exists an open cluster in $B_n$ with at least $\kappa n^d$ points.

We denote by $N_n^{(j)}$ the maximal number of open left-right crossings of $B_n$ in the $j$th direction, that is, for the maximal number of disjoint open paths connecting $B_n^{(j,-)} := \{x \in B_n : x_j = 1\}$ to $B_n^{(j,+)} := \{x \in B_n : x_j = n\}$. Recall (cf. Theorem 7.68 and Lemma 11.22 in [13]) that if $p$ is sufficiently large, then there exist positive quantities $\kappa(p)$ and $\gamma(p)$ such that for each $j$ and each $n \geq 1$, we have

$$P(N_n^{(j)} < \kappa(p) n^{d-1}) \leq e^{-\gamma(p) n^{d-1}}. \tag{B.1}$$

Note that in [13], the proof is given for bond percolation, but it can be easily adapted to site percolation. Moreover, from the proof, it is simple to derive that $\lim_{p \uparrow 1} \kappa(p) = 1$.

LEMMA B.1. *Fix $\kappa \in (0,1)$. Then, for $p < 1$, sufficiently large, there exists a positive constant $c$ such that*

$$P(\mathcal{A}_n \cap \mathcal{B}_n \cap \mathcal{C}_n(\kappa)) \geq 1 - e^{-cn} \qquad \forall n \geq 1.$$

PROOF. Let us first prove that $\mathcal{A}_n^c$, the complement of $\mathcal{A}_n$, has exponentially small probability if $p < 1$ is sufficiently large. To this end, thanks to (B.1), we may assume that there exists an open cluster $\mathcal{C}_1$ in $B_n$ intersecting both $B_n^{(1,-)}$ and $B_n^{(1,+)}$. Clearly, $d_\infty(\mathcal{C}_1) > [n/10]$. Suppose, then, that there exists another open cluster $\mathcal{C}_2$ of diameter larger than $[n/10]$. By a rather standard Peierls-like argument, we are going to show that this implies the existence of a 2-connected closed set in $B_n$ with cardinality larger than $cn$ for some constant $c > 0$ and that this latter event has exponentially small probability.

Consider the set $B_n \setminus \mathcal{C}_1$ and denote by $A$ the $B_n$-connected component of $B_n \setminus \mathcal{C}_1$ containing $\mathcal{C}_2$. We write $A_1, A_2, \ldots, A_n$ for the remaining $B_n$-connected components of $B_n \setminus \mathcal{C}_1$. Observe that $B_n \setminus A$ is $B_n$-connected. Indeed, $B_n \setminus A$ is the disjoint union of $A_1, \ldots, A_n$ and $\mathcal{C}_1$. Since $\mathcal{C}_1$ is $B_n$-connected and each $A_i$ is $B_n$-connected to $\mathcal{C}_1$, it follows that $B_n \setminus A$ is $B_n$-connected.

Define $\partial_{\text{int}} A = \{x \in A : \exists y \in B_n \setminus A, |x - y| = 1\}$. We first observe that every $x$ in $\partial_{\text{int}} A$ is closed. Indeed, any $x \in \partial_{\text{int}} A$ has a neighbor $y \in \mathcal{C}_1$ such that $|y - x| = 1$, so $\omega(x) = 1$ would imply that $x \in \mathcal{C}_1$.

Since $A$ and $B_n \setminus A$ are both $B_n$-connected, we have that $\partial_{\text{int}} A$ is 2-connected. The proof of this fact can be derived, for example, from the arguments in Appendix A of [20].

Next, we claim that

$$d_\infty(\partial_{\text{int}} A) \geq [n/10]. \tag{B.2}$$

Since $\mathcal{C}_2 \subset A$, we have $d_\infty(A) \geq d_\infty(\mathcal{C}_2) \geq [n/10]$, that is, there exist a direction $j$ and points $x, y \in A$ such that $y_j - x_j = d_\infty(A) \geq [n/10]$ where $x_j = \min\{x'_j : x' \in A\}$ and $y_j = \max\{x'_j : x' \in A\}$. If $1 < x_j$ and $y_j < n$, then we must have $x, y \in \partial_{\text{int}} A$ and (B.2) follows. Suppose, now, that $x_j = 1$, that is, $x \in B_n^{(j,-)}$ (the case $y_j = n$ is handled in the same way). If $B_n^{(j,-)} \setminus A \neq \varnothing$, then there must exist a point $z \in B_n^{(j,-)} \setminus A$ and a point $x' \in B_n^{(j,-)} \cap A$ such that $|x' - z| = 1$. In this case, $x' \in \partial_{\text{int}} A$, $x'_j = x_j = 1$ so that $d_\infty(\partial_{\text{int}} A) \geq \|y - x'\|_\infty \geq y_j - x_j$ and (B.2) follows. It remains to check the case $x_j = 1$ with $B_n^{(j,-)} \subset A$. Note that since $\mathcal{C}_1$ intersects $B_n^{(1,\pm)}$, we must have $j \neq 1$. Hence, we can exhibit points $x', y' \in A$ with $x', y' \in B_n^{(j,-)}$ and such that $x'_1 = 1$ and $y'_1 = n$. Therefore, we are in the case that we just considered above since now $x'_1 = 1$ and $B_n^{(1,-)} \setminus A \neq \varnothing$. This completes the proof of (B.2).

The above observations prove that there must exist a closed 2-connected set in $B_n$ with $\ell_\infty$-diameter at least $[n/10]$. In particular, there exists a closed 2-connected set with cardinality at least $cn$ for some constant $c > 0$. A union bound therefore gives that the probability of this event is, for $p < 1$ and $n \in \mathbb{N}$ sufficiently large, bounded above by

$$n^d \sum_{m \geq cn} e^{-\beta(p)m} e^{\alpha(d)m} \leq e^{-(c/2)\beta(p)n},$$

where we use the facts that the number of 2-connected subsets of $B_n$ with cardinality $m$ is bounded from above by $n^d e^{\alpha(d)m}$ for a $d$-dependent constant $\alpha(d)$ and that the probability that a given subset of $B_n$ with cardinality $m$ is closed is $e^{-\beta(p)m}$ with $\beta(p) \to \infty$ as $p \to 1$. In conclusion, we have shown that $\mathcal{A}_n^c$ has exponentially small probability if $p < 1$ is sufficiently large.

We turn to the events $\mathcal{B}_n$ and $\mathcal{C}_n(\kappa)$. In what follows, we take $p < 1$ sufficiently large so that $\kappa(p) \geq \kappa$. Let $\mathcal{W}_n$ be the event that $N_n^{(j)} \geq \kappa n^{d-1}$ for all $j$, $1 \leq j \leq d$. Due to (B.1) and since $d \geq 2$, we have that $P(\mathcal{W}_n^c) \leq d e^{-\gamma(p)n}$. Note that the event $\mathcal{A}_n$ implies that all of the open left-right crossings of $B_n$ must belong to the same open cluster of $B_n$. If the event $\mathcal{W}_n$ is also verified, then this common cluster has cardinality at least $\kappa n^d$ and is a crossing cluster. Hence, for $p < 1$ sufficiently large,

$$P(\mathcal{A}_n \cap \mathcal{B}_n \cap \mathcal{C}_n(\kappa)) \geq P(\mathcal{A}_n \cap \mathcal{W}_n) \geq 1 - P(\mathcal{A}_n^c) - P(\mathcal{W}_n^c) \geq 1 - e^{-cn}$$

for a suitable positive constant $c$. $\square$

REMARK 7. A set of diameter $m \in \mathbb{N}$ has less than $(2m)^d$ points. Hence, a set of at least $\kappa n^d$ points has diameter larger than $n\kappa^{1/d}/2$. In particular, taking $(1/5)^d < \kappa < 1$ in the above lemma and applying the Borel–Cantelli lemma, one obtains that if $p < 1$ is sufficiently large, then, a.s., there exists a unique maximal open cluster in $B_n$ for $n$ sufficiently large. Moreover, this unique maximal open cluster is crossing and has at least $\kappa n^d$ points.

In what follows, in order to simplify the notation, we assume that given an integer $n \geq 1$, the number $c_1 \log n$ is an integer and $n$ is a multiple of $c_1 \log n$. It is simple to adapt the result to the general case.

PROPOSITION B.2. *Given $n \geq 1$, partition $B_n$ into cubes $C_j$, $j \in J$, of side length $c_1 \log n$. Consider a maximal cluster $M(n)$ in $B_n$. The following then holds for $p < 1$ sufficiently large: there exist positive constants $\kappa$, $c'$ such that if $c_1 \geq c'$, then, a.s.,*

$$|M(n) \cap C_j| \geq \kappa |C_j| \qquad \forall j \in J,$$

*for $n$ sufficiently large.*

PROOF. Let us consider the family $\{Q_i, i \in I\}$ of cubes of side length $c_1 \log n$ included in $B_n$ given by the cubes $C_j$, $j \in J$, and by the cubes included in $B_n$ obtained my translating each $C_j$ by a distance $[c_1 \log n/2]$ in all coordinate directions. Note that $|I| \leq 2d|J| = 2dn^d/(c_1 \log n)^d$.

Fix $\kappa$ such that $(1/5)^d < \kappa < 1$. By Lemma B.1 and the Borel–Cantelli lemma, the following then holds for $p < 1$ sufficiently large and $c_1 \geq c_1^*(p)$: a.s., eventually in $n$, $B_n$ has a unique maximal open cluster $M(n)$ and, moreover, this cluster is the unique open crossing cluster and the unique open cluster of diameter larger than $[n/10]$; in addition, each cube $Q_i$ has exactly one open cluster $\mathcal{S}_i$ of diameter larger than $[c_1 \log n/10]$, and this cluster is a crossing cluster with at least $\kappa(c_1 \log n)^d$ points. Due to this characterization, it is simple to check that if $Q_i \cap Q_{i'} \neq \varnothing$, then there exists in $Q_i \cap Q_{i'}$ an open path of length larger than $[c_1 \log n/10]$ which must necessarily be included in both $\mathcal{S}_i$ and in $\mathcal{S}_{i'}$, hence $\mathcal{S}_i \cap \mathcal{S}_{i'} \neq \varnothing$. This implies that all sets $\mathcal{S}_i$ are contained in the same open cluster of $B_n$, which, moreover, has diameter $n$. Hence, this common cluster must be $M(n)$. Moreover,

$$|M(n) \cap Q_i| \geq |\mathcal{S}_i| \geq \kappa(c_1 \log n)^d \qquad \forall i \in I. \qquad \square$$

**Acknowledgments.** We thank Fabio Martinelli and Tobias Kuna for useful comments and discussions. Thanks also to an anonymous referee for helpful suggestions. We acknowledge the support of M.I.U.R. (Cofin).

Dip. Matematica
Università di Roma Tre
L. go S. Murialdo 1
00146 Roma
Italy
E-mail: caputo@mat.uniroma3.it

Dip. Matematica
Università di Roma "La Sapienza"
P. le Aldo Moro 2
00185 Roma
Italy
E-mail: faggionato@mat.uniroma1.it